\documentclass[11pt]{article}

\usepackage[T1]{fontenc}
\usepackage[margin=1in]{geometry}
\usepackage{amsmath}
\usepackage{amssymb}
\usepackage{amsthm}
\usepackage{newtxtext}
\usepackage{newtxmath}
\usepackage{bm}

\usepackage{algorithm}
\usepackage{algorithmic}
\usepackage[short]{optidef}
\usepackage{graphicx}
\usepackage{caption}
\usepackage[labelfont=sf]{subcaption}
\usepackage{booktabs}
\usepackage{multirow}
\usepackage{mathtools}
\usepackage{xcolor}
\usepackage{tikz}
\usetikzlibrary{arrows.meta}

\usepackage[sort]{natbib}
\bibpunct[, ]{(}{)}{,}{a}{}{,}

\usepackage[hidelinks]{hyperref}

\newtheorem{theorem}{Theorem}
\newtheorem{lemma}{Lemma}
\newtheorem{proposition}{Proposition}
\newtheorem{corollary}{Corollary}
\newtheorem{assumption}{Assumption}
\theoremstyle{definition}
\newtheorem{definition}{Definition}

\newcommand{\epc}{\hspace{1pc}}

\newcommand{\wt}{\widetilde}
\newcommand{\wh}{\widehat}

\DeclareMathOperator*{\argmin}{arg\,min}
\providecommand{\Bar}{\overline}

\allowdisplaybreaks
\title{Learning-Based Surrogate Method for Stochastic Optimization under
Decision-Dependent Uncertainty with Adaptive Random Designs}

\author{%
Boyang Shen \quad Junyi Liu\\[0.5em]
\normalsize Department of Industrial Engineering, Tsinghua University\\
\normalsize \href{mailto:sby22@mails.tsinghua.edu.cn}{\texttt{sby22@mails.tsinghua.edu.cn}},
\href{mailto:junyiliu@tsinghua.edu.cn}{\texttt{junyiliu@tsinghua.edu.cn}}
}
\date{}

\begin{document}

\maketitle

\begin{abstract}
We study stochastic programs in which the latent decision-dependent uncertainty is described via a nonparametric regression model. The major challenge is that, without convexity assumptions on either the cost function or the regression model, the resulting objective is both nonconvex and nonsmooth, and its first-order information is unavailable due to the unknown decision-dependent distribution. To address this issue, we construct a learning-based surrogate model that integrates simulation and statistical learning by embedding Jacobian estimates of the regression function, which are updated iteratively and interactively during the optimization procedure. We develop an adaptive random design that concentrates design points around the current iterate for Jacobian estimation and we show that the mean squared error of Jacobian estimates achieves a dimension-independent convergence rate. Building on this, we propose the learning-based stochastic prox-linear (L-SPL) algorithm with adaptive random design and establish its nonasymptotic convergence rates under various parameter settings. Numerical results demonstrate that L-SPL algorithm significantly improves sample efficiency and achieves substantially lower objective values compared to the state-of-the-art methods. More broadly, our method implies that the statistical design in an iterative learning-based optimization algorithm can be novelly tailored to the local information of the optimization procedure to sharpen estimates and enhance the convergence performance and sample efficiency of the resulting algorithm.
\end{abstract}


\section{Introduction}
\label{sec:intro}

Stochastic programming (SP) provides a broad methodology for optimizing decisions under uncertainty where uncertain quantities are modeled as random vectors \citep{birge2011introduction,shapiro2021lectures}. When uncertainty is independent of decisions, SP models can be effectively approximated and solved using discretized scenarios generated by Monte Carlo simulation.
In this paper, we study stochastic programs with \emph{decision-dependent} uncertainty (SP-DDU),  which arise in a broad range of applications including revenue management \citep{cooper2006models},  facility location \citep{basciftci2021distributionally},  appointment scheduling \citep{homem2022simulation}, and performative prediction \citep{perdomo2020performative}. Specifically, we consider the SP-DDU problem of the following form:
\begin{mini}{x \in \mathcal{X}}
    {f(x) \triangleq \mathbb{E}_{\xi \sim \mathbb{P}_\xi(x)} \left[ \varphi (x, \xi) \right] = \int \varphi(x, \xi)\, \mathbb{P}_\xi (x, \mathrm{d}\xi) ,}
    {\label{prob:intro_spddu}}
    {}
\end{mini}
where $\mathcal{X} \subseteq \mathbb{R}^d$ is a nonempty closed convex set, $\varphi$ is the cost function, and $\xi\in\mathbb{R}^\ell $ is a random vector whose probability distribution $\mathbb{P}_\xi (x)$ depends on the decision~$x$. Due to the coupled relation between the decision and the random vector, SP-DDU problems are generally nonconvex even when $\varphi$ is convex. Moreover, classical Monte Carlo methods that ignore the decision dependence could lead to systematic deterioration of solution performance \citep{cooper2006models,lee_newsvendor-type_2012}.

Without specifying a known or parametric form for the distributional map $\mathbb{P}_\xi(\,\cdot\,)$, we model the latent decision dependence through a nonparametric regression structure. Specifically, we assume that the cost function $\varphi(\,\cdot\,,\,\cdot\,)$ is weakly convex (Definition~\ref{def:weak_cvx}) and that the random vector $\xi$ satisfies
\begin{equation}\label{eq:intro_reg-model}
    \xi = c(x) + \varepsilon,
\end{equation}
where the regression function $c(x)$ is unknown and continuously differentiable, and the distribution of the random vector $\varepsilon$, denoted by $\mathbb{P}_\varepsilon$, is unknown and independent of $x$. With the regression model \eqref{eq:intro_reg-model}, \citet[Theorem 5.2]{ccinlar2011probability} implies that the problem \eqref{prob:intro_spddu} is equivalent to the following nonconvex and nonsmooth stochastic program,
\begin{mini}{x \in \mathcal{X}}
    {f(x) = \mathbb{E}_{\varepsilon\sim\mathbb{P}_\varepsilon} \left[ \varphi \left(x, c(x) + \varepsilon\right) \right] = \int \varphi(x, c(x) + \varepsilon) \mathbb{P}_\varepsilon (\mathrm{d}\varepsilon),}
    {\label{prob:intro_splddr}}
    {}
\end{mini}
in which the objective function  is weakly convex (Proposition~\ref{prop:weak_cvx_obj}) under appropriate conditions.

To solve the above problem, we propose a learning-based surrogate method by integrating simulation, learning and optimization, which is termed as the learning-based stochastic prox-linear (L-SPL) algorithm.  This algorithm iteratively constructs a prox-linear surrogate model embedded with a linear prediction model for refining the iterates along the iterative process. Different from the prominent prox-linear algorithm  \citep{davis_stochastic_2019,duchi2018stochastic}, one distinguishing feature of the L-SPL algorithm is that it incorporates simulation and learning steps in the construction of the surrogate models. At the current iterate point $z^t$, it requires the generation of decision--uncertainty data pairs $\{(X_i,\xi_i)\}$ at each iteration to obtain a statistical estimate $\wh{\nabla}c(z^t)$ of $\nabla c(z^t)$, which, together with i.i.d.\ samples $\{\eta_i\}$ from $\mathbb{P}_\xi(z^t)$, defines a linear prediction model for the local decision dependence around $z^t$. For obtaining the design points $\{X_i\}$, following the literature of nonparametric statistical learning, a standard manner is that $\{X_i\}$ follow a given distribution that is fixed during the algorithmic process, which can be referred to as the \textit{static random design}.  However, it usually suffers from the curse-of-dimensionality because the mean squared error (MSE) convergence rate of the estimator typically deteriorates with the decision dimension~$d$ \citep{stone1980optimal}.

An essential research question is whether the simulation of decision--uncertainty data pairs together with the prediction could facilitate the algorithmic process with stable and efficient convergence for solving the SP-DDU problem.  In consideration of the fact that only local information about the decision dependence around the current iterate is required by  L-SPL, our first contribution is to propose a novel \textit{adaptive random design} such that the distribution of the design points $\{X_i\}$ is adapted to this local information requirement which enables better estimation error control and the algorithm performance enhancement.
Under the adaptive random design,  the distribution of $\{X_i\}$ is centered at the iterate point $z^t$ and becomes increasingly concentrated as the sample size grows.
For the local linear regression (LLR) method \citep{fan_local_2018} under the adaptive random design, we establish a dimension-independent MSE convergence rate for the estimator $\wh{\nabla}c(z^t)$.
Building on this, the convergence analysis of the L-SPL algorithm yields a substantially improved sample complexity bound compared with the static random design.
More broadly, our results suggest that the statistical design in an iterative learning-and-optimization algorithm can be tailored to the local information requirement of the optimization procedure to sharpen estimates and improve the sample efficiency of the resulting algorithm.

Our second contribution is that we establish the nonasymptotic convergence guarantees for L-SPL to nearly stationary points in expectation (NSE, Definition~\ref{def:NSE}) of problem~\eqref{prob:intro_splddr} under both random design types.
Building on the analysis for weakly convex stochastic optimization \citep{davis_stochastic_2019,deng2021minibatch}, the major difficulty of the analysis is that, with the new scheme combining statistical estimation and surrogate approximation, we investigate the learning-induced surrogate approximation error, which involves the product of the Jacobian estimation error and the displacement from the current iterate, and we  deploy a new analytical treatment by exploiting this multiplicative structure  and the incremental sampling strategy.
We derive convergence rate of L-SPL under increasing sample sizes and tolerance-independent proximal parameters, which yields an iteration complexity of $\wt{O}(\nu^{-2})$ and a sample complexity of $\wt{O}(\nu^{-5})$ under adaptive random designs for obtaining an NSE point with the target tolerance~$\nu$. This is in contrast to the aforementioned works, which require the proximal parameter to be
a constant depending on the target tolerance~$\nu$ measuring the near-stationarity, and thus our analysis relieves the substantial tuning burden associated with the constant parameter setting.
We perform numerical experiments for two representative problems---a joint production and pricing problem and a facility location problem. The results demonstrate that L-SPL algorithm exhibits stabler convergence performance and efficient iteration and sample complexities, compared to the state-of-the-art stochastic zeroth-order method.

\textbf{Related works.} Existing algorithmic approaches for SP-DDU differ in how they handle the decision-dependent distribution. One common route is model-driven: the distributional dependence $\mathbb{P}_\xi (\,\cdot\,)$ is assumed to be known or to follow a parametric form, see, e.g., \citep{homem2022simulation,hellemo_decision-dependent_2018,bazotte2026solving,cao2024statistical}.
Such assumptions are often difficult to verify in practice, and model misspecification may degrade the quality of the resulting decisions.
Data-driven approaches, including the present work, instead learn the unknown decision dependence from decision--uncertainty pairs $\{(X_i,\xi_i)\}$ without specifying a parametric form.
Among these approaches, \citet{liu_coupled_2022} propose a derivative-free trust region method that coordinates local statistical learning with iterative optimization, and establish its asymptotic convergence to directional stationary points of the SP-DDU problem under the stringent sample size requirement.
Our nonasymptotic analysis of L-SPL establishes sample complexity guarantees in terms of statistical estimation rates, and further characterizes the effect of the adaptive random designs on the resulting sample complexity.

\citet{perdomo2020performative} and \citet{drusvyatskiy2023stochastic} study alternating schemes for SP-DDU which repeatedly optimize the decision under the distribution induced by the preceding iterate. These schemes are shown to converge to an equilibrium point~$\bar{x}$, which is defined as the solution to  the stochastic program with the fixed distribution $\mathbb{P}_{\xi}(\bar{x})$.
As indicated by \citet[Proposition 2.1]{miller2021outside}, equilibrium points may admit poor guarantees on the optimality gap for the original SP-DDU problem. Stochastic zeroth-order algorithms are recently studied for SP-DDU by   constructing gradient estimators for a smoothed approximation of $f(x)$   \citep{liu2024twotimescalederivativefreeoptimization,hikima2025zeroth,hikima_2025_zeroth_estimators,hikima2025guided,liu2026stochastic}.
 However, due to the decision-dependent distribution, controlling the variance of zeroth-order gradient estimators could be particularly challenging as evidenced by \citet{liu2024twotimescalederivativefreeoptimization}.
In contrast, L-SPL exploits the known structure of the cost function $\varphi$ and applies statistical learning only to the latent decision dependence, which could help mitigate the impact of estimation error on the optimization step.
\citet{sun2026contextual} constructs an empirical residuals-based SAA problem by embedding a global nonparametric predictor based on a given historical data set. They derive exact MIP reformulations for the SAA problems with predictors including kNN, CART and ReLU neural networks, and establish the consistency of the optimal solutions and objective values.

\textbf{Outline.} The rest of the paper is organized as follows. Section~\ref{sec:prelim} provides preliminary results on weak convexity and stationarity of the SP-DDU problem~\eqref{prob:intro_splddr}.
Section~\ref{sec:surrogate} develops the learning-based surrogate functions, establishes their surrogate properties, and proposes an adaptive random design that yields improved MSE convergence rates for the LLR estimates.
Section~\ref{sec:analysis} presents the L-SPL algorithm with the  analysis  of the nonasymptotic convergence guarantees, and derives iteration and sample complexities under multiple parameter settings.
Section~\ref{sec:nume_expt} reports numerical results on two representative application problems. We conclude the paper in Section~\ref{sec:conclusion}.

\section{Preliminaries on weak convexity and stationarity}
\label{sec:prelim}

We first introduce the basic notations used in this paper.
For any positive integer $n$, let $[n]$ represent the set $\{1,\dots, n\}$.
The notation $\|\cdot\|$ denotes the Euclidean norm for vectors and the induced operator norm for matrices.
Given a closed convex set $\mathcal{X} \subseteq \mathbb{R}^d$, let $\operatorname{dist}(z, \mathcal{X}) \triangleq \min_{x\in \mathcal{X}} \|x - z\|$ denote the distance from $z \in \mathbb{R}^d$ to $\mathcal{X}$, and $\mathcal{N}_{\mathcal{X}} (x)$ denote the normal cone to $\mathcal{X}$ at $x \in \mathcal{X}$.
For asymptotic analysis, we use the standard notations $O(\,\cdot\,)$, $\Theta(\,\cdot\,)$, $\Omega(\,\cdot\,)$, and $\wt{O}(\,\cdot\,)$, where $\wt{O}(\,\cdot\,)$ suppresses logarithmic factors.
For a locally Lipschitz continuous function $f: \mathcal{O} \to \mathbb R$ defined on an open set $\mathcal{O} \subseteq\mathbb{R}^d$, we use $\partial f(x)$ to denote the Clarke subdifferential of $f$ at $x\in\mathcal{O}$.
Given a scalar $\rho > 0$ and a closed convex set $\mathcal{X} \subset \mathcal{O}$, we use $\operatorname{prox}_{f/\rho}(x)$ to denote the proximal operator of $f/\rho$ relative to $\mathcal{X}$, namely
\[
    \operatorname{prox}_{f/\rho}(x) \triangleq \argmin_{y\in\mathcal{X}} \left\{ f(y) + \frac{\rho}{2} \|y - x\|^2 \right\}, \quad x \in \mathbb{R}^d.
\]

\renewcommand{\theenumi}{(\alph{enumi})}

Weak convexity (see the definition below) is a prominent concept that covers a broad class of nonconvex functions, such as $C^1$-smooth functions with Lipschitz gradients, and the composition of convex functions with smooth mappings \citep{duchi2018stochastic,drusvyatskiy_efficiency_2019}. Several basic characterizations of weak convexity are presented in Proposition~\ref{prop:char_weak_cvx}. Further analysis can be found in \citet[Lemma 2.1]{davis_stochastic_2019}.

\begin{definition}
\label{def:weak_cvx}
    A function $f: \mathcal{O} \rightarrow \mathbb{R}$ defined on an open convex set $\mathcal{O} \subseteq \mathbb{R}^d$ is $\lambda$-weakly convex if $f(\,\cdot\,) + \frac{\lambda}{2}\|\,\cdot\,\|^2$ is convex on $\mathcal{O}$.
\end{definition}

\begin{proposition}\label{prop:char_weak_cvx}
For a scalar $\lambda >0$ and a locally Lipschitz continuous function $f: \mathcal{O} \to \mathbb R$ on an open convex set $\mathcal{O} \subseteq \mathbb{R}^d$, the following statements are equivalent:
\begin{enumerate}
    \item $f$ is $\lambda$-weakly convex on $\mathcal{O}$.
    \item $f(\,\cdot\,) + \frac{\lambda}{2} \| \,\cdot\, - z\|^2$ is convex on $\mathcal{O}$ for every $z \in \mathcal{O}$.
    \item For any  $x,y \in \mathcal{O}$ and $g \in \partial f(y)$, the following inequality holds
    \[
        f(x) \geq f(y)+g^\top(x-y)-\frac{\lambda}{2}\|x-y\|^2.
    \]
\end{enumerate}
\end{proposition}

We impose the following assumptions on problem~\eqref{prob:intro_splddr}, covering a broad class of SP-DDU problems with weakly convex cost functions and $C^1$-smooth nonparametric regression models.

\begin{assumption}\label{assumpt:SP-DD}
    \begin{enumerate}
        \item The feasible set $\mathcal{X}$ is closed, convex, and contained in an open convex set $\mathcal O \subseteq \mathbb{R}^d$ on which the relevant functions are defined.
        \item The cost function $\varphi : \mathcal{O}  \times \mathbb{R}^{\ell} \rightarrow \mathbb{R}$ is $\lambda$-weakly convex and $\ell_\varphi$-Lipschitz continuous.
        \item The decision-dependent random vector $\xi$ follows the regression model $ \xi = c(x) + \varepsilon$, where the random vector $\varepsilon \in \mathbb{R}^\ell$ is independent of $x$ and satisfies $\mathbb{E} [\varepsilon] = 0$, $\mathbb{E}[\|\varepsilon\|^2] < \infty$.
        \item The regression function $c : \mathcal{O} \rightarrow \mathbb{R}^{\ell}$ is $\ell_c$-Lipschitz continuous and $C^1$-smooth with $L_c$-Lipschitz continuous Jacobian.
    \end{enumerate}
\end{assumption}

Recognizing that
$\varphi \left(x, c(x) + \varepsilon\right)$ is a composite random function, we show in the following proposition that the objective function $f(x)= \mathbb{E}_{\varepsilon \sim \mathbb{P}_\varepsilon} \left[ \varphi \left(x, c(x) + \varepsilon\right) \right] $ of the SP-DDU problem~\eqref{prob:intro_splddr} is weakly convex based on Proposition~\ref{prop:char_weak_cvx}. The proof is presented in Appendix~\ref{sec:proof_weak_cvx_obj}.

\begin{proposition}\label{prop:weak_cvx_obj}
    Under Assumption~\ref{assumpt:SP-DD}, the function $f(\,\cdot\,)$ in~\eqref{prob:intro_splddr} is $\left(\ell_\varphi L_c + \lambda(1+\ell_c)^2\right)$-weakly convex.
\end{proposition}

We say a point $x \in \mathcal{X}$ is stationary for the nonconvex and nonsmooth SP-DDU problem~\eqref{prob:intro_splddr} if $0 \in \partial f(x) + \mathcal{N}_\mathcal{X} (x)$.
Since our present paper addresses the sampling-based algorithm, we introduce the concept of \textit{near stationarity in expectation} for the analysis of complexity guarantees,
which is widely adopted in complexity analysis for weakly convex programs.
\begin{definition}
\label{def:NSE}
    We say a random point $x(\omega)$ is $(\nu, \delta)$-nearly stationary in expectation (NSE) for the problem~\eqref{prob:intro_splddr}, if $x\in\mathcal{X}$ holds almost surely (a.s.) and there exists a random variable $y(\omega)$ such that $y\in\mathcal{X}$ a.s., $\mathbb{E}[\|x-y\|]\leq\delta$ and $\mathbb{E} [ \operatorname{dist} (0, \partial f(y) + \mathcal{N}_\mathcal{X} (y) ) ] \leq \nu$.
\end{definition}

According to \citet{davis_stochastic_2019}, weakly convex functions naturally admit a measure for near stationarity related to Moreau envelopes, which is summarized in the proposition below.

\begin{proposition} \label{prop:stationary}
    Let $f : \mathcal{O} \rightarrow \mathbb{R}$ be a $\lambda$-weakly convex function defined on an open convex set $\mathcal{O} \subseteq \mathbb{R}^d$. Given a scalar $\rho > \lambda$ and a closed convex set $\mathcal{X} \subseteq \mathcal{O}$, define the Moreau envelope of $f$ relative to $\mathcal{X}$ by
    \begin{equation*}
        f_{1/\rho} (x) \triangleq \inf_{y \in \mathcal{X}} \left\{ f(y) + \frac{\rho}{2} \|y - x\|^2 \right\}, \quad x \in \mathbb{R}^d.
    \end{equation*}
    Then $f_{1/\rho}$ is $C^1$-smooth, and $\nabla f_{1/\rho}(x) = \rho \left(x - \operatorname{prox}_{f/\rho}(x)\right)$.
    Moreover, with $\hat{x} \triangleq \operatorname{prox}_{f/\rho}(x)$, we have $\operatorname{dist} \left(0, \partial f(\hat x) + \mathcal{N}_\mathcal{X} (\hat x) \right) \leq \rho\|\hat{x} - x\| = \|\nabla f_{1/\rho} (x)\|$.
\end{proposition}

Therefore, for a $\lambda$-weakly convex program and a candidate solution $x$,  with any scalar $\rho > \lambda$, $\left\|x - \operatorname{prox}_{f/\rho}(x)\right\|$ is commonly applied as a measure for complexity analysis \citep{davis2019proximally,davis_stochastic_2019}. Namely,   a random point $x(\omega)$
is $(\nu,\delta)$-NSE if $\mathbb E \left[\, \left\|x - \operatorname{prox}_{f/\rho}(x)\right\|\,\right] \leq \min\{\delta, \nu/\rho\}$.

\section{Surrogate models for the SP-DDU problem}
\label{sec:surrogate}

In Section~\ref{sec:concept_surrogate}, we first introduce the conceptual surrogate function based on linearizing the inner regression model, and then construct the learning-based surrogate function using statistical estimates of the regression Jacobian.
Section~\ref{sec:static_design} discusses the nonparametric statistical estimation under static random designs and illustrates the dimension dependence of the MSE convergence rates via LLR estimator.
In Section~\ref{sec:est_error_analysis}, we propose the adaptive random design and establish an improved, dimension-independent MSE convergence rate for the LLR Jacobian estimator, thereby providing tighter control of the learning-induced surrogate error.

\subsection{Learning-based surrogate functions}
\label{sec:concept_surrogate}

Given a reference point $z\in\mathcal{X}$ and a realization $\eta \sim \mathbb{P}_\xi (z)$,
we define the conceptual surrogate function $\wt{F}(\,\cdot\,,\eta;z)$ as follows,
\begin{equation}\label{eq:theo_surro}
    \wt{F} (x, \eta; z) \triangleq \varphi \left(x, \eta +\nabla c(z) (x-z)\right).
\end{equation}
We summarize several properties of the conceptual surrogate function $\wt{F} (\,\cdot\,, \eta; z)$ in the proposition below with the proof presented in Appendix~\ref{sec:proof_pro_theo_surro}.
\begin{proposition}\label{prop:pro_theo_surro}
    Under Assumption~\ref{assumpt:SP-DD}, for any reference point $z\in \mathcal X$, the surrogate function $\wt{F} (\,\cdot\,, \eta; z)$ defined in~\eqref{eq:theo_surro} satisfies the following properties:
    \begin{enumerate}
        \item Weak convexity: $\wt{F} (\,\cdot\,, \eta; z)$ is $\wt{\lambda}$-weakly convex where $\wt{\lambda} \triangleq \lambda(1+\ell_c)^2$.
        \item Quadratic surrogate error: for any $x \in \mathcal{X}$, $\left|\, \mathbb{E}\left[\,\wt{F}(x, \eta; z)\,\right] - f(x)\, \right| \leq \frac{\ell_\varphi L_c}{2}\|x-z\|^2$.
    \end{enumerate}
\end{proposition}

For the SP-DDU problem \eqref{prob:intro_splddr} with the unknown regression model, the surrogate function $\wt F(\,\cdot\,, \eta; z)$
remains conceptual. Therefore, as a natural extension, we next introduce the construction of computationally tractable surrogate functions based on the approximation of the unknown Jacobian $\nabla c(z)$ by statistical estimates.

We consider a simulated data set $S_n \triangleq \{(X_i, \xi_i)\}_{i \in [n]}$ of size $n$, where $X_i$ is a random design point and $\xi_i$ is independently generated from $\mathbb P_\xi(X_i)$ conditional on $X_i$.
We postpone the static random design and adaptive random design for generating the design points $\{X_i\}_{i \in [n]}$ in the next two subsections. Given the reference point $z$, based on the data set $S_n$, suppose that we obtain a statistical estimate $\wh{\nabla} c_n(z)$ of the true Jacobian $\nabla c(z)$.
We further generate i.i.d.\ samples $W_m \triangleq \{\eta_i\}_{i \in [m]}$ from the distribution $\mathbb{P}_\xi (z)$. Then the learning-based mini-batch surrogate function is constructed as follows:
\begin{equation}\label{eq:prac_batch_surrogate}
    \wh{F}^{n,m} \left( x,W_m; z \right) \triangleq \frac{1}{m} \sum_{i \in [m]}  \varphi \left(x, \eta_i + \wh{\nabla}c_n(z) (x-z)\right).
\end{equation}
Notice that the above surrogate is affected by the design points $\{X_i\}_{i \in [n]}$  only through the statistical estimate $\wh{\nabla}c_n(z)$.
The estimation accuracy of $\wh{\nabla}c_n(\,\cdot\,)$ can be measured via its MSE, a standard accuracy measure in statistical estimation. Since the reference point $z$ corresponds to the current iterate in the L-SPL algorithm and varies over the feasible region as the algorithm proceeds, we consider the following uniform MSE, which will be needed in the convergence analysis,
\begin{equation*}
    e_1(n) \triangleq \sup_{z\in\mathcal{X}} \sqrt{\mathbb E\left[\left\|\wh{\nabla} c_n(z) - \nabla c(z)\right\|^2\right]}.
\end{equation*}
To keep the subsequent convergence analysis applicable to different nonparametric estimators and random designs, we impose the following generic assumption that the estimation error $e_1(n)$ converges to zero at a prescribed rate with respect to the sample size $n$.

\begin{assumption}\label{assumpt:estimators}
    There exist constants $r\in(0,\frac{1}{2})$ and $C_e>0$ such that $e_1(n) \leq C_e n^{-r}$ for all $n \geq 1$.
\end{assumption}

Assumption~\ref{assumpt:estimators} is essential for controlling the learning-induced surrogate error $e_1(n)\|x - z\|$. In Assumption \ref{assumpt:estimators}, the requirement on the exponent $r<1/2$ is consistent with the common root-MSE convergence rates in nonparametric regression literature.
We examine attainable values of the convergence rate $r$ under static and adaptive random designs in the next two subsections.

We next establish the surrogate properties of $\wh F^{n,m}(\,\cdot\,,W_m;z)$ analogous to those of the conceptual surrogate function under the following boundedness assumption on the estimate $\wh{\nabla} c_n(z)$. The proof is presented in Appendix~\ref{sec:proof_prac_surro}.

\begin{assumption}\label{assu:bnd_gra_est}
    There exists a scalar $\wh{\ell}_c$ such that $\big\|\wh{\nabla} c_n(z)\big\| \leq \wh{\ell}_c$ holds almost surely for all $z\in\mathcal{X}$ and $n \geq 1$.
\end{assumption}

\begin{proposition}\label{prop:cvx_prac_surro}
    Under Assumptions~\ref{assumpt:SP-DD} and~\ref{assu:bnd_gra_est}, for any reference point $z \in \mathcal{X}$, the surrogate function $\wh{F}^{n, m} \left(\,\cdot\,,W_m; z\right)$ defined in~\eqref{eq:prac_batch_surrogate} satisfies the following properties:
    \begin{enumerate}
        \item Weak convexity: $\wh{F}^{n, m} \left(\,\cdot\,,W_m; z\right)$ is $\wh{\lambda}$-weakly convex where $\wh{\lambda} \triangleq \lambda(1+\wh{\ell}_c)^2$.
        \item Surrogate error: for any $x \in \mathcal{X}$,
        \begin{equation}\label{eq:prac_surro_error}
            \left| \mathbb{E}\left[\wh{F}^{n, m} \left(x,W_m; z\right)\right] - f(x) \right| \leq \ell_\varphi e_1(n) \|x - z\| + \frac{\ell_\varphi L_c}{2} \|x - z\|^2.
        \end{equation}
    \end{enumerate}
\end{proposition}

Notice that in the right-hand side of~\eqref{eq:prac_surro_error}, the learning-induced surrogate approximation error  involves the product of the Jacobian estimation error and the displacement from the current iterate, and such a multiplicative structure will be exploited in the convergence analysis. Moreover, any estimator satisfying Assumption~\ref{assumpt:estimators} can be truncated to satisfy Assumption~\ref{assu:bnd_gra_est} without changing its MSE convergence rate; see Lemma~\ref{lem:truncated_gradient_estimator} in Appendix~\ref{sec:estimation}. Thus in the next two subsections, we only need to address the verification of Assumption~\ref{assumpt:estimators} under static and adaptive random designs.

\subsection{Statistical estimation under static random designs}
\label{sec:static_design}

We first consider the static random design, under which the design points $\{X_i\}_{i\in [n]}$ are i.i.d.\ generated from a common density $\mathcal{M}(\,\cdot\,)$ whose support contains the feasible region $\mathcal{X}$. This is a standard random design setting in the nonparametric regression literature \citep{fan_local_2018,gu_multivariate_2015,jin_adaptive_2015}, which is generally subject to dimension-dependent convergence rates.
For a $p$-smooth regression function on $\mathbb{R}^d$, \citet{stone1980optimal} shows that the optimal rate for estimating its first-order derivatives at a point is $n^{-(p-1)/(2p+d)}$ under standard regularity conditions.
In this section, we illustrate the dimension-dependent phenomenon using local linear regression (LLR), which motivates the development of the adaptive random design. We refer to Appendix~\ref{sec:estimation} for further examples of nonparametric estimation methods.

Let $\{(X_i ,\xi_i)\}_{i \in [n]}$ be samples generated under a static random design. Given a sequence of bandwidth parameters $h_n \downarrow 0$ and an even symmetric kernel $K : \mathbb{R}^d \to \mathbb{R}_+$
(i.e., $K(x) = K(-x)$ for any $x \in \mathbb{R}^d$), at a reference point $z$, LLR obtains the estimates $\wh{c}_n(z)$ and $\wh\nabla c_n(z)$ through the following weighted least squares problem,
\begin{equation}\label{eq:llr_mean_esti}
    \left(\wh{c}_n(z), \wh\nabla c_n(z)\right) \triangleq
    \argmin_{\substack{b\in\mathbb{R}^\ell,\,A\in\mathbb{R}^{\ell\times d}}}
    \sum_{i=1}^n K \left(\frac{X_{i}-z}{h_n}\right)
    \left\| \xi_i - b - A(X_i-z) \right\|^2.
\end{equation}
The bandwidth may vary across dimensions, but we use a common bandwidth $h_n$ here for ease of exposition.
According to the asymptotic expansions for the conditional bias and covariance of the estimator $\wh\nabla c_n(z)$ \cite[Theorem~1]{lu1996multivariate},   we have
\begin{equation}\label{eq:static_conditional_mse}
    \mathbb E\left[
        \left.
        \left\|\wh\nabla c_n(z)-\nabla c(z)\right\|^2
        \,\right|\,
        X_1,\ldots,X_n
    \right]
    =
    O_p\left(
        h_n^4+\frac{1}{nh_n^{d+2}}
    \right).
\end{equation}
The above result conditioning on the design points does not necessarily imply Assumption~\ref{assumpt:estimators}.
We provide a rigorous analysis to justify Assumption~\ref{assumpt:estimators} under the proposed adaptive random design in the next subsection, by utilizing the truncation (see \eqref{eq:truncated LLR})  and well-conditioned-event decomposition techniques.
Such techniques can also be applied to estimators under static random designs for the justification of Assumption~\ref{assumpt:estimators}, while we omit the formal statement here without diverging from our goal to identify the dimension-dependent phenomenon.

Balancing the bias and variance terms in~\eqref{eq:static_conditional_mse} with $h_n=O(n^{-1/(d+6)})$ gives $r=2/(d+6)$. Under the static random design, the dimension-dependence of the variance term in \eqref{eq:static_conditional_mse}  can be explained via the effective sample size of the problem~\eqref{eq:llr_mean_esti}. In the case of the kernel $K$ with a bounded support,  design points lying in an $O(h_n)$-neighborhood of $z$ receive nonzero weights in problem~\eqref{eq:llr_mean_esti}, and the probability that a design point falls in this neighborhood is of order $h_n^d$ under a fixed density $\mathcal M$ that is continuous at $z$.
Hence, among the $n$ design points only $O(nh_n^d)$ observations in expectation  could effectively contribute to the LLR estimate $\wh\nabla c_n(z)$, which thus explains the dimension dependence by using the design points generated over the entire feasible space.

The above discussion illustrates that the performance of nonparametric estimates can be largely affected by the distribution of the design points. Combining with the fact that the L-SPL algorithm only needs the Jacobian estimates at the iterate points rather than over the entire feasible region, we propose the following adaptive random design for updating the statistical estimates adaptively along with the algorithmic process.

\subsection{Adaptive random designs and dimension-independent MSE rates}
\label{sec:est_error_analysis}

We propose the adaptive random design where the marginal distribution of the design points $\{X_i\}_{i \in [n]}$ is centered at the reference point $z$ with controllable concentration, and establish that the statistical estimates under the adaptive random design obtain the dimension-independent MSE convergence rate.
Specifically, let $\wt g : \mathbb{R}^d \rightarrow \mathbb{R}_+$ be an even symmetric probability density function. Given a reference point $z$ and sample size $n$, the adaptive random design is to generate the design points $\{X_i\}_{i\in[n]}$ following an adaptive density $\mathcal{M}(\,\cdot\,; z, n)$ by translating and scaling from $\wt g(\,\cdot\,)$ as follows:
\begin{equation}\label{eq:adaptive_design}
    \mathcal{M}(\,\cdot\,; z, n) \triangleq \frac{1}{h_n^d} \wt g \left( \frac{\,\cdot\, - z}{h_n}\right).
\end{equation}
An equivalent representation of $X_i$ is that $X_i = z + h_n U_i$, where $\{U_i\}_{i\in[n]}$ are i.i.d.\ with density~$\wt g$.
Since the weight in the LLR estimate~\eqref{eq:llr_mean_esti} satisfies $K\left(\frac{X_i-z}{h_n}\right)=K(U_i)$, in contrast to the static random design, the probability that a design point receives a nonzero weight is independent of $h_n$, which provides an intrinsic intuition of the dimension-independent MSE convergence rate under the adaptive density \eqref{eq:adaptive_design}.

Under the adaptive random design, we construct the truncated LLR estimate $\wh\nabla c_n(z)$ with a fixed scalar $\wh \ell_c\geq \ell_c$  as follows,
    \begin{equation}
    \label{eq:truncated LLR}
    \begin{array}{l}
          (\bar c_n(z),G_n(z)) \, \in \,   \underset{(b,A)\in\mathbb R^\ell\times\mathbb R^{\ell\times d}}{\mbox{argmin}}
        \sum_{i=1}^n K(U_i)
        \left\|\xi_i-b-h_nAU_i\right\|^2,\\[0.1in]
         \wh\nabla c_n(z)
        \triangleq
        \min\left\{1,\frac{\wh \ell_c}{\|G_n(z)\|}\right\}G_n(z).
        \end{array}
    \end{equation}
 In the following theorem, we provide a rigorous quantification of the dimension-independent MSE convergence rate of $\wh\nabla c_n(z)$ under the proposed adaptive random design.  In particular, the truncation guarantees a finite second moment of the estimator by enhancing its stability against rare realizations in which the least-squares problem is ill-conditioned.  The proof is provided in Appendix~\ref{sec:proof_adap}.
We expect that similar results also hold for alternative regression models, yet we postpone the rigorous verification to future works.

\begin{theorem}\label{thm:adaptive}
    Suppose that Assumption~\ref{assumpt:SP-DD} holds.
    Let $\{h_n\}_{n\geq1}$ be a positive deterministic sequence satisfying $h_n\downarrow0$, and define $\bar h\triangleq\sup_{n\geq1}h_n<\infty$.
    Given any $z\in\mathcal X$ and $n\geq1$, suppose that the following conditions are satisfied:
    \begin{enumerate}
        \item The density $\wt g$ is even symmetric and has compact support $\mathcal S \triangleq\operatorname{supp}(\wt g) \subset\mathbb R^d$.
        \item The set
            $\mathcal X^+
            \triangleq
            \left\{z+tu:z\in\mathcal X,\ u\in\mathcal S,\ t\in[0,\bar h]\right\}$
        is contained in $\mathcal O$.
        The regression function~$c$ is $C^3$-smooth on $\mathcal O$,
        and its second- and third-order derivatives are uniformly bounded there.
        \item The kernel $K:\mathbb R^d\to\mathbb R_+$ is even symmetric and bounded on $\mathcal S$, and there exists $\kappa>0$ such that $ \int_{\mathcal S}K(u)\wt g(u)\psi(u)\psi(u)^\top\,{\rm d}u \succeq\kappa\mathbb I_{d+1}$
        with $
            \psi(u)\triangleq \begin{pmatrix}1\\u\end{pmatrix} \in \mathbb{R}^{d+1}$.
    \end{enumerate}
     Let  the design points $\{X_i\}_{i\in[n]}$ follow the adaptive density $\mathcal{M}(\,\cdot\,; z, n)$ in \eqref{eq:adaptive_design}, and the truncated LLR estimates are constructed following \eqref{eq:truncated LLR}.   Then there exist positive constants $C_0$, $C_1$, and $c_0$, independent of $z$ and $n$, such that
    \begin{equation}\label{eq:adaptive_uniform_mse}
        \sup_{z\in\mathcal X}
        \mathbb E\left[
            \left\|\wh\nabla c_n(z)-\nabla c(z)\right\|^2
        \right]
        \leq
        C_0\left(h_n^4+\frac{1}{nh_n^2}\right)+C_1\exp(-c_0n).
    \end{equation}

\end{theorem}

The above theorem implies that if $h_n=h_0n^{-1/6}$ for some $h_0>0$, then there exists a constant $C_e>0$ such that
        $e_1(n)\leq C_e n^{-1/3}$ for all $n\geq1$, which justifies that Assumption~\ref{assumpt:estimators} holds with $r=1/3$.

\section{L-SPL algorithm and nonasymptotic convergence analysis}
\label{sec:analysis}

In Section~\ref{sec:algorithm}, we develop the learning-based stochastic prox-linear (L-SPL) algorithm for solving the SP-DDU problem~\eqref{prob:intro_splddr}, based on the learning-based surrogate function constructed in Section~\ref{sec:surrogate}.
Section~\ref{sec:pre_analysis} establishes a nonasymptotic convergence guarantee to $(\nu,\delta)$-NSE points.
We then derive the iteration and sample complexities under both constant and variable parameter settings in Section~\ref{sec:conv_guarantee}.

\subsection{L-SPL algorithm}
\label{sec:algorithm}

We present the L-SPL algorithm in Algorithm \ref{alg:L-SPL} with the process described as follows.
Each iteration of the main loop consists of three phases: simulation, learning, and optimization. At iteration $t$, in the simulation phase, we generate $m_t$ i.i.d.\ samples from $\mathbb{P}_\xi(z^t)$, which is called the mini-batch data set denoted by $W^t_{m_t}$. We  also generate the data set $S^t_{n_t}$ of decision--uncertainty data pairs of size $n_t$ with the design points following a specified random design. The random design is either a static random design that is fixed throughout the algorithm, or the adaptive random design that updates adaptive density $\mathcal M(\cdot;z^t,n_t)$ along the iterative process.  In the learning phase, we construct the Jacobian estimate $\wh \nabla c^t_{n_t}(z^t)$  based on $S^t_{n_t}$. The optimization phase constructs the learning-based surrogate function according to~\eqref{eq:prac_batch_surrogate} and obtains the minimum solution to the corresponding proximal subproblem utilizing a convex programming solver. After completing $T+1$ iterations, the algorithm generates the output solution by sampling from the sequence of iteration points according to a specified discrete distribution.

\begin{algorithm}[ht]
    \caption{Learning-based Stochastic Prox-Linear (L-SPL) Algorithm}
    \label{alg:L-SPL}
    \begin{algorithmic}
        \REQUIRE $z^0 \in \mathcal{X}$, $T$, random design type $\mathcal{D}\in\{\mathrm{static},\mathrm{adaptive}\}$, $\rho > \ell_\varphi L_c + \max\{\wt{\lambda}, \wh{\lambda}\}$, $\{\alpha_t\}_{t=0}^T \subset (\rho + \ell_\varphi L_c, +\infty)$, $\{m_t\}_{t=0}^T$, $\{n_t\}_{t=0}^T$.

        \FOR{$t = 0, \cdots, T$}
        \STATE \textsc{Simulation phase:} Generate $W^t_{m_t} \triangleq \{\eta_i^t\}_{i\in[m_t]}$ from $\mathbb{P}_\xi(z^t)$.
        \IF{$\mathcal{D}=\mathrm{adaptive}$}
            \STATE Generate $S^t_{n_t} \triangleq \{(X_i^t,\xi_i^t)\}_{i\in[n_t]}$ with $X_i^t$ drawn from $\mathcal{M}(\,\cdot\,;z^t,n_t)$ and $\xi_i^t$ drawn from $\mathbb{P}_\xi(X_i^t)$.
        \ELSE
            \STATE Generate $S^t_{n_t} \triangleq \{(X_i^t,\xi_i^t)\}_{i\in[n_t]}$ with $X_i^t$ drawn from $\mathcal{M}(\,\cdot\,)$ and $\xi_i^t$ drawn from $\mathbb{P}_\xi(X_i^t)$.
        \ENDIF

        \STATE \textsc{Learning phase:} Estimate $\wh{\nabla}c^t_{n_t}(z^t)$ based on $S^t_{n_t}$.

        \STATE \textsc{Optimization phase:} Construct $\wh{F}^{n_t,m_t} \left( \,\cdot\,, W^t_{m_t}; z^t \right)$ according to~\eqref{eq:prac_batch_surrogate}.
        Compute
        \[
            z^{t+1} = \argmin_{z\in\mathcal{X}} \left\{\wh{F}^{n_t,m_t} \left(z, W^t_{m_t}; z^t\right) + \frac{\alpha_t}{2}\|z - z^t\|^2\right\}.
        \]
        \ENDFOR

        \ENSURE $z^{t^*}$ with $t^* \in \{0,\dots,T\}$ sampled according to $\{p_t\}_{t=0}^{T}$ defined by $p_t \propto \frac{1}{\alpha_t + \rho - \ell_\varphi L_c - \wt{\lambda} - \wh{\lambda}}$.
    \end{algorithmic}
\end{algorithm}

The L-SPL algorithm is governed by three groups of algorithmic parameters: the sample sizes $\{n_t\}_{t=0}^T$, $\{m_t\}_{t=0}^T$, and the proximal parameters $\{\alpha_t\}_{t=0}^T$. These parameters control the statistical accuracy, mini-batch stability, and proximal regularization, respectively. We next study the nonasymptotic convergence guarantee of the L-SPL algorithm.

\subsection{Nonasymptotic convergence analysis}
\label{sec:pre_analysis}

Recall from Proposition~\ref{prop:stationary} that, with $\wt{\nu} \triangleq \min\{\delta, \nu/\rho\}$ and the proximal point $\wh{z}^{\,t} \triangleq {\rm prox}_{f/\rho}(z^t)$ for a parameter $\rho > \ell_\varphi L_c + \lambda(1 + \ell_c)^2$, the following condition is sufficient for $z^{t^*}$ to be a $(\nu, \delta)$-NSE point:
\begin{equation}\label{eq:suff_NSE_cond_vari}
    \mathbb{E} \, \left[ \, \left\| z^{t^*} - \wh{z}^{\,t^*} \right\| \, \right] \leq \wt{\nu}.
\end{equation}
Let $B^t \triangleq S^t_{n_t} \cup W^t_{m_t}$ denote the samples used at iteration $t$, and let $\mathbb{E}_{B^t} [\,\cdot\,]$ denote conditional expectation given all past samples $B^0, \dots, B^{t-1}$.
In Proposition \ref{prop:opt_error}, we aim to quantify the optimization error in expectation $\left| \, \mathbb{E}_{B^t} \left[\, \wh{F}^{n_t,m_t} \left(z^{t+1}, W_{m_t}^t; z^t\right) - f(z^{t+1})\, \right] \,  \right|$. To address the intrinsic dependence of the new iterate $z^{t+1}$ on samples $B^t$, we decompose the optimization error into the stability estimate and the stochastic surrogate error, and further obtain the bound by utilizing Proposition~\ref{prop:cvx_prac_surro} and the mini-batch stability estimate of \citet[Theorem~3.2]{deng2021minibatch}.
The proof is presented in Appendix~\ref{sec:proof_opt_error}.
\begin{proposition}\label{prop:opt_error}
    Under Assumptions~\ref{assumpt:SP-DD} and \ref{assu:bnd_gra_est}, suppose that $\alpha_t > \lambda(1 + \wh{\ell}_c)^2$.
    Then we have
    \begin{equation*}
    \begin{aligned}
        &\left| \, \mathbb{E}_{B^t} \left[\, \wh{F}^{n_t,m_t} \left(z^{t+1}, W_{m_t}^t; z^t\right) - f(z^{t+1})\, \right] \,  \right| \\
        &\qquad \leq \frac{2\ell_\varphi^2(1+\wh{\ell}_c)^2}{m_t(\alpha_t - \lambda(1 + \wh{\ell}_c)^2)} + \ell_\varphi e_1(n_t) \sqrt{\mathbb{E}_{B^t}\left[ \left\|z^{t+1} - z^t\right\|^2 \right]} + \frac{\ell_\varphi L_c}{2}\mathbb{E}_{B^t}\left[ \left\|z^{t+1} - z^t\right\|^2 \right] .
    \end{aligned}
    \end{equation*}
\end{proposition}

For the surrogate function $\wh{F}^{n_t,m_t} \left(\,\cdot\,, W_{m_t}^t; z^t\right)$, the learning-induced surrogate-error terms (e.g., the second term in the right-hand-side of the above inequality) in Propositions~\ref{prop:cvx_prac_surro} and~\ref{prop:opt_error}
distinguish our analysis from those considered by \citet{davis_stochastic_2019} and \citet{deng2021minibatch}.
The next proposition establishes a descent property for the Moreau envelope $f_{1/\rho}(\,\cdot\,)$ while accounting for the coupling between the estimation error $e_1(n_t)$ and the terms $\|\wh z^{\,t}-z^t\|$ and $\sqrt{\mathbb{E}_{B^t}[\|z^{t+1}-z^t\|^2]}$.
In particular, the term $e_1(n_t)\|\wh z^{\,t}-z^t\|$ coupling the estimation error with the near stationarity measure  is retained on the left-hand side of~\eqref{eq:des_property} and thus requires new treatments for establishing the convergence guarantee in Theorem~\ref{thm:nonasym_basic}.
The proof is deferred to Appendix~\ref{sec:proof_des_property}.

\begin{proposition}
\label{prop:des_property}
    Under Assumptions~\ref{assumpt:SP-DD} and \ref{assu:bnd_gra_est}, denote $\tau \triangleq \ell_\varphi L_c + \wh{\lambda} + \wt{\lambda}$. Suppose that $\rho > \ell_\varphi L_c + \max\{\wt{\lambda},\wh{\lambda}\}$ and $\{\alpha_t\}_{t=0}^T \subset (\rho + \ell_\varphi L_c, +\infty)$.
    Then the sequence of iteration points $\left\{ z^t \right\}_{t=0}^T$ produced by Algorithm~\ref{alg:L-SPL} satisfies
    \begin{equation}\label{eq:des_property}
        \begin{aligned}
        &\frac{\rho}{\alpha_t + \rho - \tau} \mathbb{E}\left[\frac{2\rho - \ell_\varphi L_c - \tau}{2} \left\|\, \wh{z}^{\,t} - z^t \,\right\|^2 - \ell_\varphi e_1(n_t) \left\|\, \wh{z}^{\,t} - z^t \,\right\|\right] \\
        &\qquad\leq \mathbb{E} \left[ f_{1/\rho} \left(z^t\right) - f_{1/\rho} \left(z^{t+1}\right) \right] + \frac{\rho}{\alpha_t + \rho - \tau} \left[\frac{2 \ell_\varphi^2 (1 + \wh{\ell}_c)^2}{m_t(\alpha_t - \wh{\lambda})} + \frac{ \ell_\varphi^2 e_1^2(n_t)}{2(\alpha_t - \rho - \ell_\varphi L_c)}\right].
        \end{aligned}
    \end{equation}
\end{proposition}

We next obtain the fundamental nonasymptotic guarantee under variable settings of the parameters $\{\alpha_t\}$, $\{m_t\}$ and $\{n_t\}$, with the proof provided in Appendix~\ref{sec:proof_nonasym_basic}.

\begin{theorem}\label{thm:nonasym_basic}
    Under Assumptions~\ref{assumpt:SP-DD} and \ref{assu:bnd_gra_est}, denote $\Delta_0 \triangleq f_{1/\rho} \left(z^0\right) - \min_{z\in\mathcal{X}} f(z)$. Suppose that $\rho > \ell_\varphi L_c + \max\{\wt{\lambda},\wh{\lambda}\}$ and $\{\alpha_t\}_{t=0}^T \subset (\rho + \ell_\varphi L_c, +\infty)$.
    Let $t^* \in \left\{0, \dots, T\right\}$ be a random variable following the discrete distribution $\{p_t\}_{t=0}^T$ with $p_t \propto \frac{1}{\alpha_t + \rho - \tau}$.
    Then the solution $z^{t^*}$ returned by Algorithm~\ref{alg:L-SPL} satisfies
    \begin{equation*}
        \mathbb{E}\left\|\, \wh{z}^{\,t^*} - z^{t^*} \,\right\| \leq \sqrt{\frac{4\Delta_0 + \sum_{t=0}^T \frac{\rho}{\alpha_t + \rho - \tau} \left[\frac{8 \ell_\varphi^2 (1 + \wh{\ell}_c)^2}{m_t(\alpha_t - \wh{\lambda})} + \frac{4 \ell_\varphi^2 e_1^2(n_t)}{2\rho - \ell_\varphi L_c - \tau} + \frac{2 \ell_\varphi^2 e_1^2(n_t)}{\alpha_t - \rho - \ell_\varphi L_c}\right] }{\sum_{t=0}^T \frac{\rho(2\rho - \ell_\varphi L_c - \tau)}{\alpha_t + \rho - \tau}}} .
    \end{equation*}
\end{theorem}

\subsection{Iteration and sample complexities}
\label{sec:conv_guarantee}

In this section, we derive iteration and sample complexities for achieving $\mathbb{E}\left\|\wh{z}^{\,t^*} - z^{t^*}\right\| = O(\wt{\nu})$ with the target tolerance $\wt{\nu}=\min\{\delta,\nu/\rho\}$ as stated in~\eqref{eq:suff_NSE_cond_vari}. The complexities are derived under constant and variable parameter settings, distinguished by whether their specification depends on the target tolerance.
We first analyze parameters that remain constant across iterations, requiring values to be tolerance-dependent.
We then consider parameters increasing at the polynomial rate such that they are independent of the target tolerance.
The resulting complexity bounds illustrate how the statistical estimation rate affects the performance of L-SPL.

\subsubsection{Complexity under constant parameter settings}
\label{sec:limited_vari_param}

We first consider the constant parameter setting with $m_t \equiv m > 0$, $n_t \equiv n > 0$, and $\alpha_t \equiv \alpha > \rho+\ell_\varphi L_c$ as a baseline. By Theorem~\ref{thm:nonasym_basic}, we have
\begin{equation}\label{eq:const_bnd}
    \mathbb{E}\left\|\, \wh{z}^{\,t^*} - z^{t^*} \,\right\| \leq \sqrt{\frac{\frac{4\Delta_0(\alpha + \rho - \tau)}{\rho(T+1)} + \frac{8\ell_\varphi^2 (1 + \wh{\ell}_c)^2}{m(\alpha - \wh{\lambda})} + \frac{4 \ell_\varphi^2 e_1^2(n)}{2\rho - \ell_\varphi L_c - \tau} + \frac{2 \ell_\varphi^2 e_1^2(n)}{\alpha - \rho - \ell_\varphi L_c}}{2\rho - \ell_\varphi L_c - \tau} }.
\end{equation}
Under Assumption~\ref{assumpt:estimators}, with $r$ denoting the MSE convergence rate, the right-hand side of~\eqref{eq:const_bnd} is in the order of $O\left(\sqrt{T^{-1}\alpha + m^{-1}\alpha^{-1} + n^{-2r} + n^{-2r} \alpha^{-1}}\right)$.
Thus the iteration complexity is  $O(\alpha\wt{\nu}^{-2})$, provided that $n = \Theta \left(\wt{\nu}^{-1/r}\right)$ and $m\alpha = \Theta \left(\wt{\nu}^{-2}\right)$, which are all tolerance-dependent.  The latter requirement on  $m\alpha$ leads to two natural constant parameter settings. If $m=O(1)$ and $\alpha=\Theta(\wt{\nu}^{-2})$, which is analogous to the constant parameter settings in \citet{deng2021minibatch} and \citet{davis_stochastic_2019}, then we obtain the iteration complexity of $O(\wt{\nu}^{-4})$ and the sample complexity of $O(\wt{\nu}^{-4-1/r})$. Alternatively, if $m=\Theta(\wt{\nu}^{-2})$ and $\alpha=O(1)$, then the iteration complexity and the sample complexity reduce to $O(\wt{\nu}^{-2})$ and $O(\wt{\nu}^{-2-1/r})$ by the fact that $r<1/2$.
However, the dependence on the target tolerance $\wt{\nu}$ makes them difficult to be specified in practice and creates a substantial tuning burden.

\subsubsection{Complexity under variable parameter settings}
\label{sec:best_vari_param}

We next consider variable parameter settings that progressively reduce the stochastic errors without depending on the target tolerance~$\wt{\nu}$. Specifically, let
$\alpha_t = \alpha_0 (t+1)^b$, $m_t = \lceil m_0(t+1)^j \rceil$, and $n_t = \lceil n_0(t+1)^k \rceil$, where $\alpha_0 > \rho + \ell_\varphi L_c$, $m_0, n_0 \in \mathbb{Z}_+$, and $b,j,k \geq 0$.
As indicated by Theorem~\ref{thm:nonasym_basic}, a larger growth rate~$b$ of the proximal parameter damps stochastic errors but weakens the cumulative descent governed by $\sum_{t=0}^T \frac{\rho(2\rho - \ell_\varphi L_c - \tau)}{\alpha_t + \rho -   \tau}$.
In the following corollary, by fixing $b$, we identify the choices of the sample size growth rates~$j$ and~$k$ with the tightest complexity. The proof is presented in Appendix~\ref{sec:proof_coroll_vari}.

\begin{corollary}\label{coroll:vari}
    With $\rho > \ell_\varphi L_c + \max\{\wh{\lambda}, \wt{\lambda}\}$, let $\alpha_t = \alpha_0 (t+1)^b$, $m_t = \lceil m_0(t+1)^j \rceil$, $n_t = \lceil n_0(t+1)^k \rceil$ where $\alpha_0 > \rho+\ell_\varphi L_c$, $m_0, n_0 \in \mathbb{Z}_+$ and $b, j, k \geq 0$.
    Then under Assumptions~\ref{assumpt:SP-DD}--\ref{assu:bnd_gra_est}, for any $b \in [0, \frac{1}{2}]$, by minimizing the bound derived from Theorem~\ref{thm:nonasym_basic}, the choices
    \begin{equation*}
        j^*(b) \triangleq 1 - 2b, \qquad k^*(b) \triangleq \frac{1-b}{2r}
    \end{equation*}
    attain the tightest sample complexity $\wt{O}\left( \wt{\nu}^{-\frac{2}{1-b}-\frac{1}{r}} \right)$ and the corresponding  iteration complexity $\wt{O}\left( \wt{\nu}^{-\frac{2}{1-b}} \right)$  for obtaining a $(\nu,\delta)$-NSE solution $z^{t^*}$ by Algorithm~\ref{alg:L-SPL}.
\end{corollary}

Table~\ref{tab:comp_const_param} summarizes the complexity for the variable settings obtained from Corollary~\ref{coroll:vari} with $b=\frac{1}{2}$ and $b=0$, together with the complexity results for the constant parameter settings in Section~\ref{sec:limited_vari_param}.
Corollary~\ref{coroll:vari} indicates that within the polynomial family of parameter settings, setting $b=0$ leads to the lowest iteration and sample complexities by theory.
Though the theoretical guarantees can be conservative,  they suggest a parameter selection rule: keep the proximal parameter a sufficiently large constant and incrementally increase sample sizes for the proper control of the stochastic errors and also convergence efficiency. This is also evidenced in numerical experiments in Section \ref{sec:numeric_pp}.

\begin{table}[htbp]
\centering
\caption{Complexities for obtaining a $(\nu, \delta)$-NSE point by Algorithm~\ref{alg:L-SPL}.}
\label{tab:comp_const_param}
\renewcommand{\arraystretch}{1.2}
\resizebox{\textwidth}{!}{%
\begin{tabular}{cccccc}
    \toprule
    Parameter setting & $\{n_t\}$ & $\{m_t\}$ & $\{\alpha_t\}$ & Iteration complexity & Sample complexity \\
    \midrule

    Constant I & $\Theta\left(\wt{\nu}^{-1/r}\right)$ & $O(1)$ & $\Theta\left(\wt{\nu}^{-2}\right)$ & $O\left(\wt{\nu}^{-4}\right)$ & $O\left(\wt{\nu}^{-4-1/r}\right)$ \\

    \addlinespace[0.4em]

    Constant II & $\Theta\left(\wt{\nu}^{-1/r}\right)$
    & $\Theta\left(\wt{\nu}^{-2}\right)$ & $O(1)$
    & $O\left(\wt{\nu}^{-2}\right)$ & $O\left(\wt{\nu}^{-2-1/r}\right)$ \\

    \addlinespace[0.4em]

    Variable I ($b=\frac{1}{2}$) & $n_0(t+1)^{1/(4r)}$ & $m_0$ & $\alpha_0(t+1)^{1/2}$ & $\wt{O}\left(\wt{\nu}^{-4}\right)$ & $\wt{O}\left(\wt{\nu}^{-4-1/r}\right)$ \\

    \addlinespace[0.4em]

    Variable II ($b=0$) & $n_0(t+1)^{1/(2r)}$ & $m_0(t+1)$ & $\alpha_0$ & $\wt{O}\left(\wt{\nu}^{-2}\right)$ & $\wt{O}\left(\wt{\nu}^{-2-1/r}\right)$ \\

    \bottomrule
\end{tabular}
}
\smallskip
\begin{minipage}{\textwidth}
\footnotesize
\textit{Note.} $\wt{\nu} = \min\{\delta, \nu/\rho\}$ with $\rho > \ell_\varphi L_c + \lambda(1 + \ell_c)^2$, and $r$ denotes the convergence rate in Assumption~\ref{assumpt:estimators}.
For LLR under the adaptive random design in Theorem~\ref{thm:adaptive}, one can take $r=1/3$; under static random designs, the attainable rate $r$ typically deteriorates with the dimension $d$.
\end{minipage}
\end{table}

  In Table~\ref{tab:comp_const_param}, the variable settings I and II which are tolerance independent  match the complexity orders of the constant settings I and II, respectively, while allowing the parameters to be selected independently of $\wt{\nu}$.
Moreover, the choice between static and adaptive random designs affects the sample complexity through the estimation rate~$r$. In particular, Theorem~\ref{thm:adaptive} gives $r=1/3$ for LLR under the specified adaptive random design, whereas the rate under static random designs typically suffers dimensional degradation.

\section{Numerical results}
\label{sec:nume_expt}

We present numerical results on two representative problems, a joint production and pricing problem and a facility location problem, in Sections~\ref{sec:numeric_pp} and~\ref{sec:numeric_facility}, respectively.
For the first problem, we examine the performance of L-SPL under multiple parameter settings and random design types.
For the second problem, we compare L-SPL with the stochastic zeroth-order (SZO) method proposed by \citet{hikima_2025_zeroth_estimators} under the same sample budgets.

Unless otherwise stated, we implement L-SPL using the adaptive random design~\eqref{eq:adaptive_design}, where $\wt g$ is the uniform density over $[-1,1]^d$.
We use LLR with the kernel $K(x) = \left(3/4\right)^d \prod_{j\in[d]} \max \left\{1-x_j^2, 0\right\}$ to estimate the Jacobian $\nabla c(z^t)$.
The numerical experiments were conducted on a Windows 11 Pro desktop with an Intel Core i5-14600KF processor (3.50 GHz) and 32 GB RAM. The implementation uses Python 3.11.13, with Gurobi 13.0.0 used to solve the convex subproblems in the L-SPL algorithm.
The code is available at \url{https://github.com/Shenby18/spddu-lspl}.

\subsection{Joint production and pricing problem}
\label{sec:numeric_pp}
We consider a joint production and pricing problem with two products. For each product $i=1,2$, the decision maker determines the production quantity $q_i$ and unit price $p_i$. Let $p \triangleq (p_1,p_2)\in\mathbb{R}_+^2$ and $q \triangleq (q_1,q_2)\in\mathbb{R}_+^2$. The random demand $D(p)\in\mathbb{R}_+^2$ depends on the price through the regression model $D(p)=\phi(p)+\varepsilon$, where the components of $\varepsilon\in\mathbb{R}^2$ are independent and uniformly distributed over $[-\Bar{\varepsilon}_i,\Bar{\varepsilon}_i]$, respectively. The regression function $\phi(\,\cdot\,)=(\phi_1(\,\cdot\,),\phi_2(\,\cdot\,))$ is given by
\begin{equation*}
    \phi_i(p) = \frac{u_i \cdot \exp{(v_i - w_i p_i)}}{1 + \exp{(v_1 - w_1 p_1)} + \exp{(v_2 - w_2 p_2)}} + \Bar{\varepsilon}_i, \quad i = 1, 2,
\end{equation*}
where $u_i$, $v_i$, and $w_i$ are positive scalars, and the additive term $\Bar{\varepsilon}_i$ ensures that the demand is nonnegative. Let $c_1\in\mathbb{R}_+^2$ denote the regular production cost per unit, $c_2\in\mathbb{R}_+^2$ the expedited production cost per unit for unmet demand, and $c_3\in\mathbb{R}_+^2$ the holding cost per unit for leftover products. The joint production and pricing problem is formulated as
\[
\begin{aligned}
    \min_{p, q} \quad &c_1^\top q + \mathbb{E}_{D(p) \sim \mathbb{P}_D(p)} \left[ -p^\top D(p) + c_2^\top \max\{D(p) - q, 0\} + c_3^\top \max\{q - D(p), 0\} \right] \\
    \text{s.t.} \quad &0 \leq q \leq q_{\rm max}, \quad 0 \leq p \leq p_{\rm max}.
\end{aligned}
\]
The problem parameter values are summarized in Table~\ref{tab:jpp_prob_param}.

\begin{table}[htbp]
\centering
\caption{Parameter values for the joint production and pricing problem.}
\label{tab:jpp_prob_param}
\begin{tabular}{l@{\qquad}ccccccccc}
        \toprule
        Parameter&  $c_1$&  $c_2$&  $c_3$&  $u$&  $v$&  $w$&  $\Bar{\varepsilon}$&  $p_{\rm max}$&  $q_{\rm max}$\\
        \midrule
        Value&  $(3,2)$&  $(7.5,9)$&  $(3,3)$&  $(6,10)$&  $(7,8)$&  $(1,0.8)$&  $(1,1)$&  $(10,10)$&  $(15,15)$\\
        \bottomrule
\end{tabular}
\end{table}

To investigate the computational performance of L-SPL, we specify different algorithmic settings, each of which is tested over 50 independent replications. Within each replication, the initial point is generated uniformly from the feasible region and shared across settings.
We measure performance by the relative optimality gap $|f(z^t)-f^*|/|f^*|$, where $f(z^t)$ can be evaluated analytically in this problem and the optimal value $f^*$ is approximated by grid search over the feasible region.
Regarding the LLR bandwidth parameter, we set $h_t\equiv h$ under the constant parameter settings and $h_t = h_0(t+1)^{-k/6}$ under the variable parameter settings with $n_t = n_0(t+1)^k$, as suggested by Theorem~\ref{thm:adaptive}. We tune the constants $h$ and $h_0$ from $\{0.5,1.0,\ldots,6.0\}$ with the best mean performance at the sample budget of $60{,}000$ in preliminary runs.
In Figures \ref{fig:const_dimi_stepsize} and \ref{fig:jpp_oracle_schedule}, the solid lines represent the median relative optimality gap over 50 replications, smoothed by a centered rolling median with window size 7 to make the trends clearer, while the shaded regions report the interquartile range of the raw results.

\subsubsection{Comparison of constant and variable parameter settings}
\label{sec:expt_stepsize}

\begin{figure}[htbp]
\centering
    \subcaptionbox{Different constant sample sizes.}[0.45\textwidth]{\includegraphics[width=0.45\textwidth]{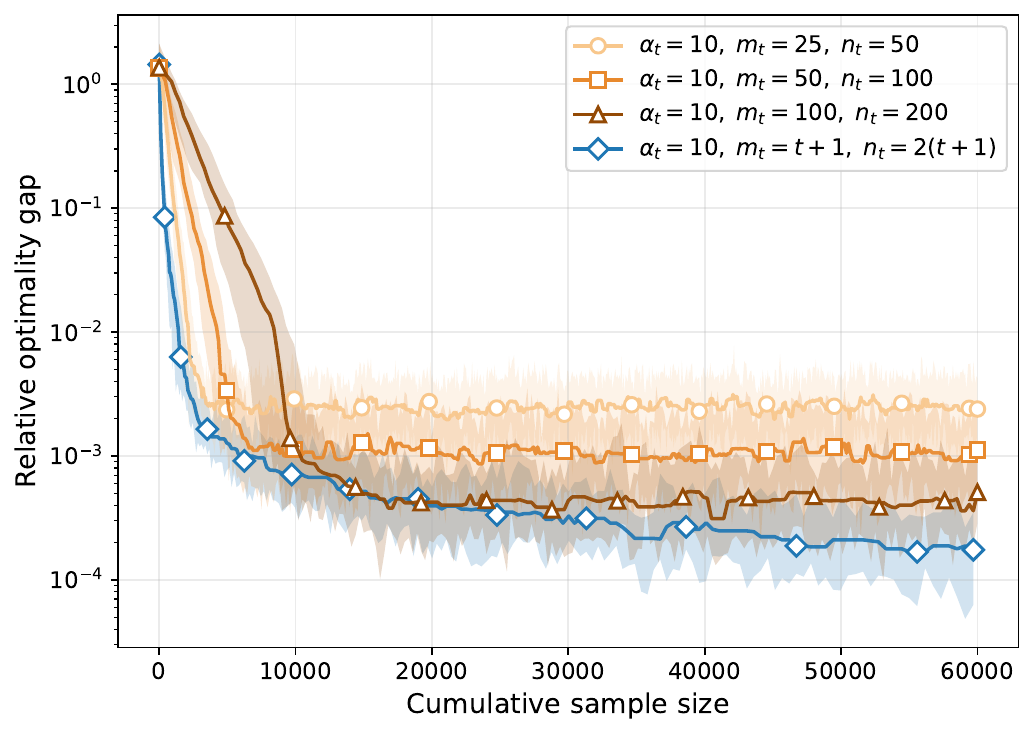}}
    \hfill
    \subcaptionbox{Different constant proximal parameters.}[0.45\textwidth]{\includegraphics[width=0.45\textwidth]{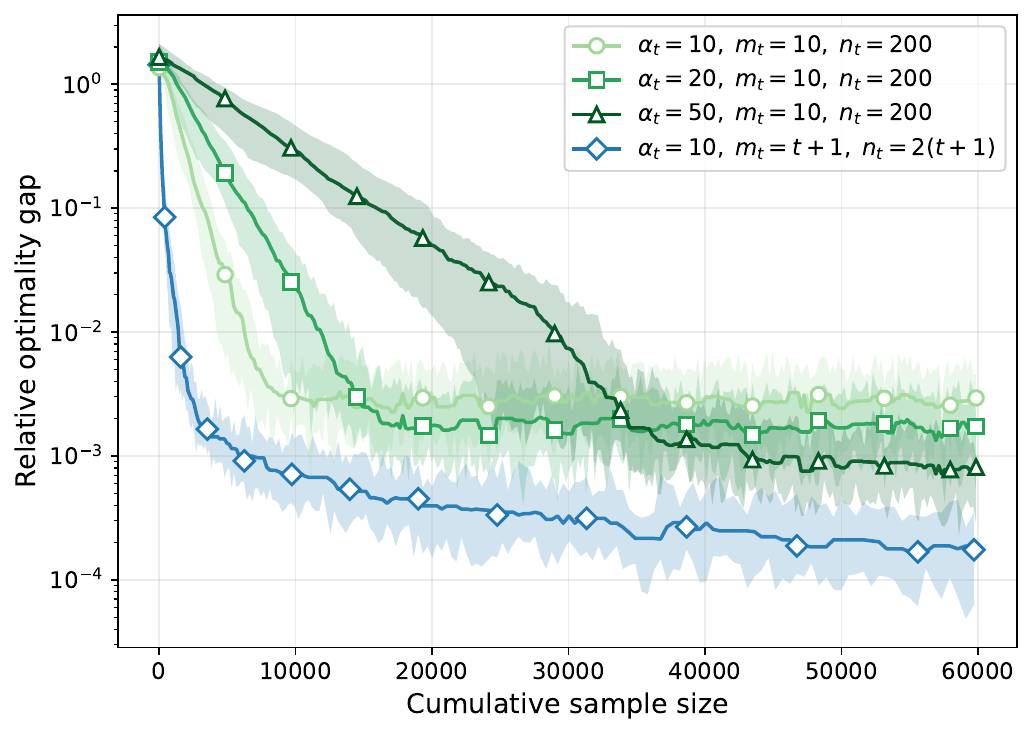}}
\caption{Convergence comparison of L-SPL with constant and variable parameter settings.}
\label{fig:const_dimi_stepsize}
\end{figure}

We first compare the convergence behavior of L-SPL under constant and variable parameter settings.
Figure~\ref{fig:const_dimi_stepsize}(a) compares different constant sample sizes while fixing $\alpha_t\equiv 10$, whereas Figure~\ref{fig:const_dimi_stepsize}(b) compares different constant proximal parameters while fixing $m_t\equiv 10$ and $n_t\equiv 200$.
Both panels include the variable setting $\alpha_t\equiv 10$, $m_t=t+1$, and $n_t=2(t+1)$ as a common reference.

Both panels show that L-SPL with constant parameter settings quickly reaches the stage dominated by the parameter-dependent error without further progress.
In contrast, the variable setting (blue curve in Figure \ref{fig:const_dimi_stepsize}) progressively reduces the stochastic errors.
Across both panels, larger constant sample sizes or proximal parameters produce lower optimality gaps but also slower improvement, particularly at early iterations, illustrating a trade-off between early sample efficiency and attainable accuracy associated with the choice of constant sample sizes and proximal parameter values.

\subsubsection{Comparison across variable parameter settings and random designs}

\begin{figure}[htbp]
\centering
    \subcaptionbox{Comparison by iteration number}{\includegraphics[width=.45\textwidth]{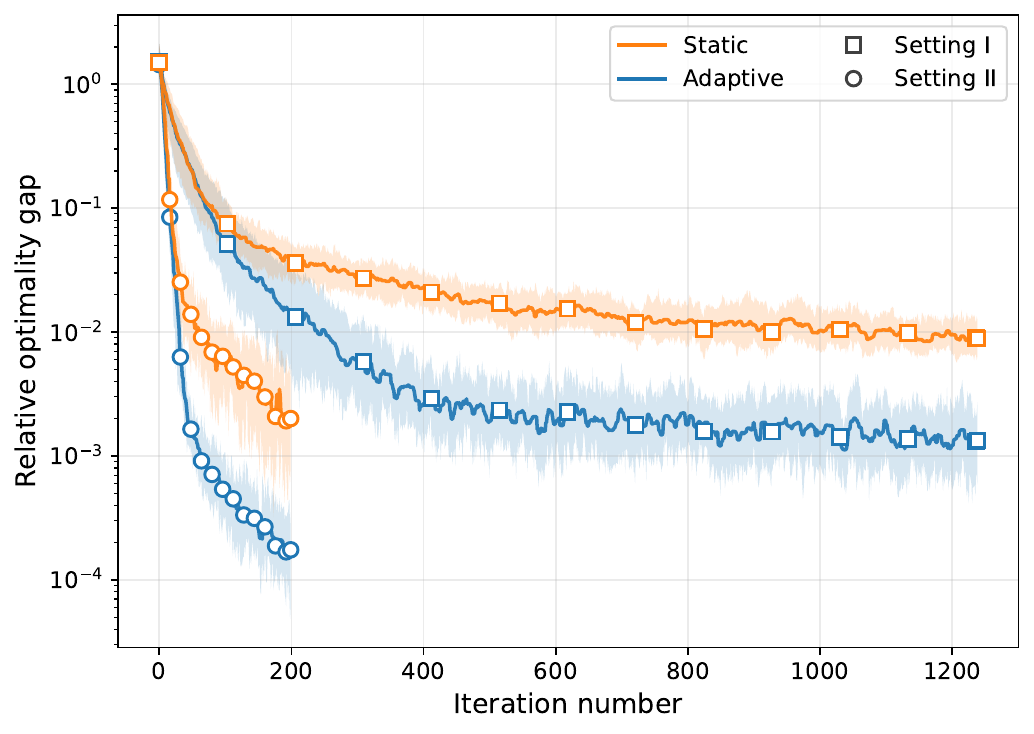}}
    \hfill
    \subcaptionbox{Comparison by cumulative sample size}{\includegraphics[width=.45\textwidth]{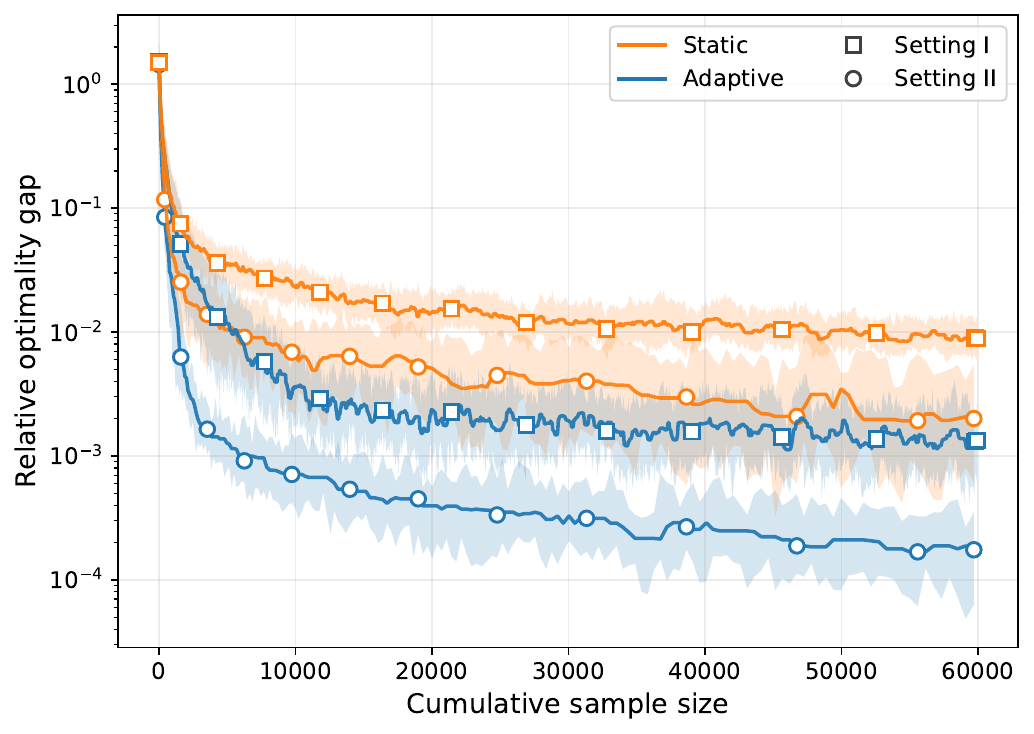}}
\caption{Convergence comparison of L-SPL under two variable parameter settings and two random designs.}
\label{fig:jpp_oracle_schedule}
\end{figure}

We next compare the two variable parameter settings in Table~\ref{tab:comp_const_param} under the two random design types. Specifically, we consider
\[
\begin{aligned}
    \text{Setting I:}\quad&\alpha_t=10(t+1)^{0.5},\qquad m_t=1,\qquad n_t=\left\lceil 2(t+1)^{0.5}\right\rceil,\\
    \text{Setting II:}\quad&\alpha_t=10,\qquad m_t=t+1,\qquad n_t=2(t+1).
\end{aligned}
\]
Both settings are tested under the adaptive random design specified before and a static random design in which the design points are generated uniformly over the feasible region.
Figure~\ref{fig:jpp_oracle_schedule} reports the relative optimality gap against the iteration number and cumulative sample size, respectively.
Each run of L-SPL is terminated when the cumulative sample usage $\sum_{s=0}^t(m_s+n_s)$ reaches the common budget of $60{,}000$.
Since Setting II uses faster-growing sample sizes, it completes fewer iterations within the common sample budget.

Figure~\ref{fig:jpp_oracle_schedule} shows that under both random design types, Setting II outperforms Setting I by attaining smaller optimality gaps within the same number of iterations and the same cumulative sample budget.
This comparison supports the practical guideline suggested in Section~\ref{sec:best_vari_param}: keeping the proximal parameter fixed and increasing the sample sizes can be more effective than increasing the proximal parameter with slowly growing sample sizes.
For either parameter setting, L-SPL under the adaptive random design attains smaller optimality gaps than under the static random design, illustrating that concentrating the design points around the current iterate improves local statistical estimation and leads to significantly better objective performance.

\subsection{Facility location problem}
\label{sec:numeric_facility}

\subsubsection{Problem setup}

We consider a facility location problem in which a service provider chooses the locations of facilities to serve customer sites.
Suppose there are $I$ customer sites, where site $i\in[I]$ is located at $(a_i,b_i)$, and the decision maker determines the locations of $J$ facilities denoted by $\{(x_j,y_j)\}_{j \in [J]}$, which affect the random demand $D_{ij} (x, y) \in \mathbb{R}_+$ from customer site $i$ for service at facility $j$ where $x \triangleq (x_j)_{j\in[J]}$ and $y \triangleq (y_j)_{j\in[J]}$.  Customers must physically reach a facility, so the commuting distance affects their willingness to receive service, which leads to the decision-dependent demand model as follows,
\begin{equation*}
    D_{ij} (x, y) = \beta_{ij} (x, y)\mu_i(x, y) + \varepsilon_{ij},
\end{equation*}
with
\[\begin{gathered}
    \beta_{ij}(x,y)  \triangleq \frac{\exp(-l_i(x_j,y_j) / \gamma)}{\sum_{j=1}^J \exp(-l_i(x_j,y_j) / \gamma)}, \quad \mu_i(x,y) \triangleq \wt{\mu}_i \left( 1-\exp\left(-\frac{1}{J}\sum_{j=1}^J \frac{1}{  u \cdot  l_i (x_j,y_j)} \right) \right), \\
    l_i(x_j,y_j) \triangleq  (x_j-a_i)^2+(y_j-b_i)^2.
\end{gathered}\]
Here $\wt{\mu}_i$, $u$, and $\gamma$ are positive scalars, and $\{\varepsilon_{ij}\}_{i \in [I], j \in [J]}$ are independent random variables uniformly distributed over $[0,\overline{\varepsilon}]$.
The term $\mu_i(x,y)$ represents the total potential demand from customer site $i$, which increases as facilities move closer to this site and saturates at $\wt{\mu}_i$.
This demand then distributes among the facilities according to the weights $\{\beta_{ij}(x,y)\}_{j\in[J]}$, with closer facilities attracting larger shares.
The increasing and saturating features of $\mu_i(\,\cdot\,, \,\cdot\,)$ follow the moment function structure of
\citet{basciftci2021distributionally}, and the form of $\beta_{ij}(\,\cdot\,, \,\cdot\,)$ is adopted from \citet{aydin2022multi}.

Let $r_i$ and $p_i$ denote the unit revenue for satisfied demand
and the unit penalty for unsatisfied demand at site $i$, respectively. Given a demand realization $D(x, y) \triangleq \{D_{ij}(x,y)\}_{i \in [I], j \in [J]}$, the maximum revenue
is determined by
\[
\begin{aligned}
    \varphi(D(x, y)) \triangleq \min_{\{z_{ij}:\, i\in [I], j \in [J]\}} \quad & \sum_{i=1}^I p_i \sum_{j=1}^J (D_{ij} (x, y) - z_{ij}) - \sum_{i=1}^I r_i \sum_{j=1}^J z_{ij} \\
    \text{s.t.} \quad& \sum_{i=1}^I z_{ij}\leq C_j, \quad \forall j \in [J] , \\
    & 0 \leq z_{ij} \leq D_{ij}(x, y), \quad \forall  i\in [I], j \in [J],
\end{aligned}
\]
where $C_j$ is the capacity of facility $j \in [J]$, and $z_{ij}$ is the served demand from site $i$ by facility $j$. The facility locations are then determined by the following stochastic program
\[
\begin{aligned}
    \min_{x, y} \quad & \mathbb{E}_{D(x, y) \sim \mathbb{P}_D (x,y)} [\, \varphi(D(x, y)) \,] \\
    \text{s.t.} \quad & \underline{x} \leq x_j \leq \overline{x}, \quad  \underline{y} \leq y_j \leq \overline{y}, \quad \forall j \in [J].
\end{aligned}
\]
In this experiment, the site locations $\{(a_i,b_i)\}_{i \in [I]}$ are sampled uniformly from the region $[\underline{x}, \overline{x}]\times[\underline{y}, \overline{y}]$. The parameters $C_j$, $\wt\mu_i$, $r_i$ and $p_i$ are sampled from uniform distributions $U(20, 40)$, $U(15, 25)$, $U(0.5, 1.5)$ and $U(1,3)$, respectively. Other problem parameter values are summarized in Table~\ref{tab:faci_prob_param}.

\begin{table}[htbp]
\centering
\caption{Parameter values for the facility location problem.}
\label{tab:faci_prob_param}
\makebox[\textwidth][c]{\begin{tabular*}{0.5\textwidth}{@{\extracolsep{\fill}}lccccccc}
        \toprule
        Parameter&  $\underline{x}$&  $\overline{x}$&  $\underline{y}$&  $\overline{y}$& $\overline{\varepsilon}$& $u$ &$\gamma$ \\
        \midrule
        Value&  0&  10&  0&  10& 2&  0.3 &5 \\
        \bottomrule
\end{tabular*}}
\end{table}

\subsubsection{Experimental results}

We compare L-SPL with SZO.
For L-SPL,
we use constant parameters $\alpha$, $m$, $n$ and $h$, which are selected via grid search with $\alpha \in \{0.1, 1, 10\}$, $m \in \{1, 5, 10\}$, $n \in \{10, 20, 30, 50\}$, and $h \in \{1, 1.5, 2, 2.5, 3\}$.
For SZO, we implement Algorithm~1 of \citet{hikima_2025_zeroth_estimators} using the gradient estimator in their equation~(28).
We tune the smoothing parameter~$\mu$, the step size $\beta$, the number of search directions $N_{\rm dir}$, and the number of samples per direction $m_{\rm sam}$ via grid search with $\mu \in \{0.05, 0.1, 0.25, 0.5, 1, 1.5\}$, $\beta \in \{0.001, 0.01, 0.1, 1\}$, $N_{\rm dir} \in \{1, 10, 20, 50, 100\}$, and $m_{\rm sam} \in \{1, 10, 20, 50, 100\}$.
We select the above parameters with the lowest mean objective at the sample budget of $5{,}000$ in preliminary runs.

Table~\ref{tab:facility_als_szo} reports the objective values
of L-SPL and SZO across different problem sizes $(I,J)$ and cumulative sample budgets $N$.
Under each problem size, we conduct 50 independent replications for one problem instance. For each replication, both methods are initialized from a collection of 10 starting points that are uniformly sampled from the feasible region, and the trajectory with the lowest validation objective at $N=5000$ is selected. This procedure is designed to mitigate the effect of highly nonconvex curvature of the objective. For each selected trajectory, we evaluate the objective value at each iterate where the cumulative sample usage first reaches the budget $N \in \{100, 500, 1000, 2000, 5000\}$, using an evaluation sample set of size 500.

\begin{table}[htbp]
\centering
\caption{Objective value
comparison between L-SPL and SZO
for the facility location problem.}
\label{tab:facility_als_szo}
\renewcommand{\arraystretch}{1.2}
\resizebox{\textwidth}{!}{%
\begin{tabular}{clccccc}
\toprule
$(I,J)$ & Metric & $N=100$ & $N=500$ & $N=1000$ & $N=2000$ & $N=5000$ \\
\midrule
\multirow{3}{*}{$(5,2)$}
& L-SPL Obj.
& $-44.64 \pm 13.80$
& $-60.88 \pm 13.86$
& $-68.59 \pm 6.92$
& $-71.29 \pm 2.40$
& $-71.69 \pm 1.99$ \\
& SZO Obj.
& $-41.91 \pm 11.64$
& $-49.67 \pm 12.22$
& $-54.35 \pm 12.36$
& $-62.20 \pm 9.92$
& $-71.16 \pm 2.15$ \\
& $\Delta_{\mathrm{rel}}$
& $\mathbf{11.039\%}$
& $\mathbf{28.755\%}$
& $\mathbf{33.040\%}$
& $\mathbf{17.990\%}$
& $\mathbf{0.780\%}$ \\
\midrule
\multirow{3}{*}{$(10,6)$}
& L-SPL Obj.
& $-148.49 \pm 17.15$
& $-175.62 \pm 17.41$
& $-187.16 \pm 13.47$
& $-198.84 \pm 8.90$
& $-205.62 \pm 3.44$ \\
& SZO Obj.
& $-139.96 \pm 16.82$
& $-162.18 \pm 16.82$
& $-178.27 \pm 15.29$
& $-191.99 \pm 10.03$
& $-202.38 \pm 3.77$ \\
& $\Delta_{\mathrm{rel}}$
& $\mathbf{7.340\%}$
& $\mathbf{9.582\%}$
& $\mathbf{5.819\%}$
& $\mathbf{3.854\%}$
& $\mathbf{1.630\%}$ \\
\midrule
\multirow{3}{*}{$(20,14)$}
& L-SPL Obj.
& $-309.15 \pm 14.49$
& $-343.53 \pm 7.51$
& $-353.71 \pm 7.50$
& $-361.07 \pm 4.72$
& $-369.09 \pm 2.36$ \\
& SZO Obj.
& $-318.97 \pm 15.77$
& $-346.02 \pm 8.41$
& $-355.42 \pm 6.48$
& $-360.69 \pm 5.19$
& $-369.22 \pm 1.93$ \\
& $\Delta_{\mathrm{rel}}$
& $-2.869\%$
& $-0.664\%$
& $-0.451\%$
& $\mathbf{0.125\%}$
& $-0.034\%$ \\
\bottomrule
\end{tabular}}
\smallskip
\begin{minipage}{\textwidth}
\footnotesize
\textit{Note.} We report the mean $\pm$ standard deviation of objective values.
$\Delta_{\mathrm{rel}}$ reports the mean of $(\mathrm{Obj}_{\mathrm{SZO}}-\mathrm{Obj}_{\mathrm{L\text{-}SPL}}) / |\mathrm{Obj}_{\mathrm{SZO}}| \times 100\%$. Thus, a positive $\Delta_{\mathrm{rel}}$ indicates that L-SPL attains a lower objective value than SZO.
\end{minipage}
\end{table}

Table~\ref{tab:facility_als_szo} suggests that L-SPL is more sample-efficient than SZO  on the small-sized and medium-sized instances,  with roughly  5\%--30\% relative improvement. The two methods perform comparably on the largest instance with $(I,J) = (20, 14)$.
One possible explanation lies in the different estimation targets of the two methods. Although the adaptive random design yields a dimension-independent MSE convergence rate, the finite-sample error of the Jacobian estimator for the $IJ$-dimensional demand map may still increase with $I$ and $J$.
In contrast, SZO constructs a zeroth-order gradient estimate of the scalar objective function.
We conjecture that the heavier learning task embedded in L-SPL for large-scale problems may offset part of the benefit of exploiting the known cost structure in the surrogate subproblems.

\section{Conclusion}\label{sec:conclusion}

We develop the L-SPL algorithm for SP-DDU problems with latent decision dependence, which integrates simulation, statistical learning and proximal surrogate optimization.
We propose an adaptive random design embedded in the L-SPL algorithm that concentrates the design points around the current iterate to improve local statistical estimation.
Using LLR as a concrete nonparametric estimation method, we establish a dimension-independent MSE convergence rate for the Jacobian estimate under the adaptive random design, whereas the corresponding rate under static random designs generally deteriorates with the decision dimension.
We establish nonasymptotic convergence guarantees for L-SPL and derive its iteration and sample complexities under constant and variable parameter settings.
Our complexity analysis shows that keeping the proximal parameter constant while increasing the sample sizes yields more favorable iteration and sample complexities than increasing the proximal parameter with slowly growing sample sizes.
The numerical experiments illustrate the advantage of L-SPL algorithm with sample efficiency and faster convergence with adaptive random design and variable parameter settings.
This work illustrates an important aspect of simulation and optimization integration for SP-DDU problems. Future research may extend the proposed framework to stochastic programs with both contextual and decision-dependent uncertainty, as well as to multistage stochastic programming and chance-constrained programming with decision-dependent uncertainty.

\appendix
\numberwithin{equation}{section}
\numberwithin{theorem}{section}
\numberwithin{lemma}{section}
\numberwithin{proposition}{section}
\numberwithin{corollary}{section}
\section{Omitted proofs}
\label{sec:proofs}

In this appendix we present the proofs omitted from the main body.

\subsection{Proof of Proposition~\ref{prop:weak_cvx_obj}}
\label{sec:proof_weak_cvx_obj}

Denote $F (x, \varepsilon) \triangleq \varphi (x, c(x) + \varepsilon)$. For any $x, y \in \mathcal{O}$ and $v \in \partial_x F(x, \varepsilon)$, the chain rule of Clarke subdifferentials implies that there exists $(w_1, w_2)\in \partial \varphi(x, c(x) + \varepsilon)$ such that $v = w_1 + \nabla c(x)^\top w_2$.
By the weak convexity of $\varphi(\,\cdot\,,\,\cdot\,)$, we obtain
\begin{equation*}
\begin{aligned}
    F(y, \varepsilon) &= \varphi(y, c(y)+\varepsilon) \geq \varphi(x, c(x)+\varepsilon) + w_1^\top (y-x) + w_2^\top (c(y)-c(x)) - \frac{\lambda}{2}\left\|\begin{matrix}
        y - x \\ c(y) - c(x)
    \end{matrix}\right\|^2 \\
    &= F(x, \varepsilon) + v^\top(y-x) + w_2^\top(c(y)-c(x)-\nabla c(x)(y-x)) - \frac{\lambda}{2}\left\|\begin{matrix}
        y - x \\ c(y) - c(x)
    \end{matrix}\right\|^2.
\end{aligned}
\end{equation*}
The Lipschitz continuity of $c(\,\cdot\,)$, $\varphi(\,\cdot\,,\,\cdot\,)$ and $\nabla c(\,\cdot\,)$ yields $\left\|\begin{matrix}
        y - x \\ c(y) - c(x)
    \end{matrix}\right\|^2 \leq (1+\ell_c)^2 \|y - x\|^2$ and
\begin{equation*}
\begin{aligned}
    w_2^\top(c(y)-c(x)-\nabla c(x)(y-x)) &\geq -\|w_2\| \|c(y)-c(x)-\nabla c(x)(y-x)\| \geq -\frac{\ell_\varphi L_c}{2}\|y-x\|^2.
\end{aligned}
\end{equation*}
Combining the above inequalities, we obtain \[
    F(y, \varepsilon) \geq F(x, \varepsilon) + v^\top(y-x) - \frac{\ell_\varphi L_c + \lambda (1+\ell_c)^2}{2}\|y-x\|^2 ,
\] which implies that $F(\,\cdot\,, \varepsilon)$ is $\left(\ell_\varphi L_c + \lambda (1+\ell_c)^2\right)$-weakly convex by Proposition~\ref{prop:char_weak_cvx}.
We complete the proof by noticing that $f(x) = \mathbb{E}[F(x,\varepsilon)]$.

\subsection{Proof of Proposition~\ref{prop:pro_theo_surro}}
\label{sec:proof_pro_theo_surro}

We first prove part (a). For any $x, y \in \mathcal{O}$ and $v \in \partial_x \wt{F}(x,\eta;z)$, the chain rule of Clarke subdifferentials implies that there exists $(w_1, w_2) \in \partial \varphi(x, \eta + \nabla c(z)(x-z))$ such that $v = w_1 + \nabla c(z)^\top w_2$. By the weak convexity of $\varphi(\,\cdot\,,\,\cdot\,)$, we have
\begin{equation*}
\begin{aligned}
    \wt{F}(y,\eta;z) &= \varphi \left(y, \eta + \nabla c(z) (y-z)\right) \\
    &\geq \varphi \left(x, \eta + \nabla c(z) (x-z)\right) + w_1^\top (y-x) + w_2^\top \nabla c(z) (y-x) - \frac{\lambda}{2}\left(1 + \|\nabla c(z)\|\right)^2\|y - x\|^2 \\
    &\geq \wt{F}(x,\eta;z) + v^\top (y - x) - \frac{\lambda}{2}(1+\ell_c)^2 \|y - x\|^2,
\end{aligned}
\end{equation*}
where the last inequality is obtained by Assumption~\ref{assumpt:SP-DD}(d). Thus $\wt{F}(\,\cdot\,,\eta;z)$ is $\lambda(1+\ell_c)^2$-weakly convex.
We next prove part (b). Since $\eta\sim\mathbb P_\xi(z)$, under the regression model in Assumption~\ref{assumpt:SP-DD}(c), there exists a realization $\varepsilon \sim \mathbb{P}_\varepsilon$ such that $\eta = c(z) + \varepsilon$. By the Lipschitz continuity of $\varphi(\,\cdot\,,\,\cdot\,)$ and $\nabla c(\,\cdot\,)$, we obtain
\begin{equation*}
\begin{aligned}
    \left|\,\mathbb E\left[\wt F(x,\eta;z)\right]-f(x)\,\right|
    &= \Big|\,\mathbb E\left[
        \varphi\left(x,c(z)+\varepsilon+\nabla c(z)(x-z)\right)
        -\varphi\left(x,c(x)+\varepsilon\right)
    \right]\,\Big| \\
    &\leq \ell_\varphi
    \left\|c(x)-c(z)-\nabla c(z)(x-z)\right\| \leq \frac{\ell_\varphi L_c}{2}\|x-z\|^2,
\end{aligned}
\end{equation*}
which completes the proof.

\subsection{Proof of Proposition~\ref{prop:cvx_prac_surro}}
\label{sec:proof_prac_surro}

The proof of part~(a) is similar to that of Proposition~\ref{prop:pro_theo_surro}(a) and is therefore omitted. For part~(b),
we deduce
\begin{equation*}
\begin{aligned}
    &\mathrel{\phantom{=}}\left|\mathbb E\left[\wh F^{n,m}(x,W_m;z)\right]-f(x)\right| \\
    &= \left|\mathbb E\left[\varphi\left(x,\eta_1+\wh{\nabla}c_n(z)(x-z)\right) - \varphi\left(x,\eta_1+\nabla c(z)(x-z)\right) + \varphi\left(x,\eta_1+\nabla c(z)(x-z)\right)\right]-f(x)\right| \\
    &\leq \left|\mathbb E\left[
        \varphi\left(x,\eta_1+\wh{\nabla}c_n(z)(x-z)\right)
        -\varphi\left(x,\eta_1+\nabla c(z)(x-z)\right)
    \right]\right| + \left|\mathbb E\left[\wt F(x,\eta_1;z)\right]-f(x)\right| \\
    &\leq \ell_\varphi e_1(n)\|x-z\|
        +\frac{\ell_\varphi L_c}{2}\|x-z\|^2,
\end{aligned}
\end{equation*}
where we obtain the first equality since the samples in $W_m$ are identically distributed, and the second inequality follows from the Lipschitz continuity of $\varphi(\,\cdot\,,\,\cdot\,)$, Proposition~\ref{prop:pro_theo_surro}(b), and the Cauchy-Schwarz inequality.

\subsection{Proof of Theorem~\ref{thm:adaptive}}\label{sec:proof_adap}

Notice that an equivalent representation of $X_i$ is that $X_i=z+h_nU_i$, where $\{U_i\}_{i\in[n]}$ are i.i.d.\ with density $\wt g$. Fix $n\geq1$ and $z\in\mathcal X$.
We first define the good event
\[
    \mathcal G_n
    \triangleq
    \left\{\|M_n-M\|\leq\frac{\kappa}{2}\right\},
\]
where
\[
    M \triangleq \int_{\mathcal S}K(u)\wt g(u)\psi(u)\psi(u)^\top\,{\rm d}u, \quad M_n\triangleq
    \frac{1}{n}\sum_{i=1}^n
    K(U_i)\psi(U_i)\psi(U_i)^\top.
\]
We have $\mathbb E[M_n]=M$.
Since $\mathcal S$ is compact and $K$ is bounded on $\mathcal S$, the matrices $K(U_i)\psi(U_i)\psi(U_i)^\top$ are uniformly bounded.
Applying Hoeffding's inequality,
there exist constants $C_M,c_M>0$ such that
\begin{equation}\label{eq:adaptive_good_event}
    \mathbb P(\mathcal G_n^c)\leq C_M\exp(-c_Mn).
\end{equation}
These constants do not depend on $z$ or $n$.
It follows from $M \succeq \kappa \mathbb{I}_{d+1}$ that $\lambda_{\min}(M)\geq\kappa$. Thus on $\mathcal G_n$, Weyl's inequality gives
    $\lambda_{\min}(M_n)
    \geq
    \lambda_{\min}(M)-\|M_n-M\|
    \geq\frac{\kappa}{2}$, and hence,
\begin{equation}\label{eq:adaptive_inverse_bound}
    \|M_n^{-1}\|\leq\frac{2}{\kappa}.
\end{equation}

We next analyze the untruncated LLR estimate $G_n(z)$ on $\mathcal G_n$.
Let
\[
    B(z)\triangleq
    \begin{pmatrix}c(z)&h_n\nabla c(z)\end{pmatrix},
    \qquad
    \bar B_n(z)\triangleq
    \begin{pmatrix}\bar c_n(z)&h_nG_n(z)\end{pmatrix},
\]
which are $\ell\times(d+1)$ matrices. The normal equations for the untruncated LLR estimate $\bar B_n(z)$ are
\[
    \frac{1}{n}\sum_{i=1}^n
    K(U_i)
    \left(\xi_i-\bar B_n(z)\psi(U_i)\right)
    \psi(U_i)^\top
    =0.
\]
With $R_z(u)\triangleq c(z+h_nu)-c(z)-h_n\nabla c(z)u$, the response $\xi_i$ can be written as
\[
    \xi_i = B(z)\psi(U_i)+R_z(U_i)+\varepsilon_i.
\]
Substituting this expression for $\xi_i$ into the normal equations gives
\[
\begin{aligned}
    \left(\bar B_n(z)-B(z)\right)M_n
    &={}
    \frac{1}{n}\sum_{i=1}^n
    K(U_i)
    \left(R_z(U_i)+\varepsilon_i\right)
    \psi(U_i)^\top = Q_n+N_n,
\end{aligned}
\]
where
\[
    Q_n\triangleq
    \frac{1}{n}\sum_{i=1}^n
    K(U_i)R_z(U_i)\psi(U_i)^\top, \qquad N_n\triangleq
    \frac{1}{n}\sum_{i=1}^n
    K(U_i)\varepsilon_i\psi(U_i)^\top.
\]
On $\mathcal G_n$, the matrix $M_n$ is nonsingular.
Let $J\triangleq\begin{pmatrix}\mathbf 0_d^\top\\\mathbb I_d\end{pmatrix}\in\mathbb R^{(d+1)\times d}$.
We have $\left(\bar B_n(z)-B(z)\right)J=h_n\left(G_n(z)-\nabla c(z)\right)$.
Therefore,
\begin{equation}\label{eq:adaptive_gradient_decomposition}
    G_n(z)-\nabla c(z)
    =
    \frac{1}{h_n}\left(N_n+Q_n\right)M_n^{-1}J
    \qquad\text{on }\mathcal G_n.
\end{equation}

We first control the noise term $\frac{1}{h_n}N_nM_n^{-1}J$.
Denote $\sigma^2\triangleq\mathbb E\|\varepsilon\|^2<\infty$.
Conditional on $\{U_i\}_{i\in[n]}$, the independence and zero mean of $\{\varepsilon_i\}_{i\in[n]}$ imply
\[
    \mathbb E\left[
        \left.\|N_n\|_{\rm F}^2\,\right|\,U_1,\ldots,U_n
    \right]
    =
    \frac{\sigma^2}{n^2}
    \sum_{i=1}^n K^2(U_i)\|\psi(U_i)\|^2,
\]
Here, $\|\cdot\|_{\rm F}$ denotes the Frobenius norm.
Define $K_{\max}\triangleq\sup_{u\in\mathcal S}K(u)$ and $R_{\mathcal S}\triangleq\sup_{u\in\mathcal S}\|u\|$.
Since $U_i\in\mathcal S$ almost surely and $\|\psi(U_i)\|^2=1+\|U_i\|^2$, we have $K^2(U_i)\|\psi(U_i)\|^2 \leq K_{\max}^2\left(1+R_{\mathcal S}^2\right)$ almost surely.
Thus, with $C_N\triangleq\sigma^2K_{\max}^2\left(1+R_{\mathcal S}^2\right)$,
the tower property and the preceding conditional second-moment imply
\[
\begin{aligned}
    \mathbb E\left[
        \|N_n\|_{\rm F}^2
    \right]
    &=
    \mathbb E\left[
        \mathbb E\left[
            \left.\|N_n\|_{\rm F}^2\,\right|\,
            U_1,\ldots,U_n
        \right]
    \right] \leq
    \frac{C_N}{n}.
\end{aligned}
\]
Since $J^\top J=\mathbb I_d$, we have $\|J\|=1$.
Hence, on $\mathcal G_n$, the inequality $\|ABC\|_{\rm F}\leq\|A\|_{\rm F}\|B\|\|C\|$ and~\eqref{eq:adaptive_inverse_bound} give
\begin{equation}\label{eq:adaptive_noise_bound}
    \mathbb E\left[
        \left\|
            \frac{1}{h_n}N_nM_n^{-1}J
        \right\|_{\rm F}^2
        \mathbf 1_{\mathcal G_n}
    \right] \leq
    \frac{1}{h_n^2} \mathbb{E}\left[\|N_n\|_{\rm F}^2
    \|M_n^{-1}\|^2\|J\|^2 \mathbf 1_{\mathcal G_n}\right]
    \leq\frac{4C_N}{\kappa^2nh_n^2}.
\end{equation}

We next control the Taylor remainder term $\frac{1}{h_n}Q_nM_n^{-1}J$.
With $Q\triangleq\mathbb E[Q_n] = \int_{\mathcal S} K(u) \wt{g}(u) R_z(u) \psi(u)^\top \mathrm{d}u$, we have
\begin{equation}\label{eq:adaptive_remain_decom}
    \frac{1}{h_n}Q_nM_n^{-1}J
    =
    \frac{1}{h_n}(Q_n-Q)M_n^{-1}J
    +\frac{1}{h_n}Q(M_n^{-1}-M^{-1})J
    +\frac{1}{h_n}QM^{-1}J.
\end{equation}
We use the following bounds, which will be established below, with $C>0$ denoting a generic constant independent of $z$ and $n$ whose value may change from line to line:
\begin{equation}\label{eq:adaptive_remainder_ingredients}
\begin{aligned}
    \mathbb E\left[\|Q_n-Q\|_{\rm F}^2\right]
    &\leq\frac{Ch_n^4}{n},
    &
    \|Q\|_{\rm F}
    &\leq Ch_n^2,\\
    \|QM^{-1}J\|_{\rm F}
    &\leq Ch_n^3,
    &
    \mathbb E\left[\|M_n-M\|^2\right]
    &\leq\frac{C}{n},\\
    \|M_n^{-1}-M^{-1}\|
    &\leq C\|M_n-M\|
    &&\text{on }\mathcal G_n.
\end{aligned}
\end{equation}
Together with~\eqref{eq:adaptive_inverse_bound} and $\|J\|=1$, these bounds give
\[
\begin{aligned}
    \mathbb E\left[
        \left\|
            \frac{1}{h_n}(Q_n-Q)M_n^{-1}J
        \right\|_{\rm F}^2
        \mathbf 1_{\mathcal G_n}
    \right]
    &\leq
    \frac{1}{h_n^2}
    \mathbb E\left[
        \|Q_n-Q\|_{\rm F}^2
        \|M_n^{-1}\|^2
        \|J\|^2
        \mathbf 1_{\mathcal G_n}
    \right] \leq
    \frac{4}{\kappa^2h_n^2}
    \mathbb E\left[\|Q_n-Q\|_{\rm F}^2\right]
    \leq\frac{Ch_n^2}{n},\\
    \mathbb E\left[
        \left\|
            \frac{1}{h_n}Q(M_n^{-1}-M^{-1})J
        \right\|_{\rm F}^2
        \mathbf 1_{\mathcal G_n}
    \right]
    &\leq
    \frac{1}{h_n^2}\|Q\|_{\rm F}^2
    \mathbb E\left[
        \|M_n^{-1}-M^{-1}\|^2
        \mathbf 1_{\mathcal G_n}
    \right] \leq
    \frac{C}{h_n^2}\|Q\|_{\rm F}^2
    \mathbb E\left[\|M_n-M\|^2\right]
    \leq\frac{Ch_n^2}{n},\\
    \left\|
        \frac{1}{h_n}QM^{-1}J
    \right\|_{\rm F}^2
    &\leq
    \frac{1}{h_n^2}\|QM^{-1}J\|_{\rm F}^2
    \leq Ch_n^4.
\end{aligned}
\]
Applying $\|A+B+C\|_{\rm F}^2\leq3(\|A\|_{\rm F}^2+\|B\|_{\rm F}^2+\|C\|_{\rm F}^2)$ to~\eqref{eq:adaptive_remain_decom}, we obtain
\begin{equation}\label{eq:adaptive_remainder_bound}
    \mathbb E\left[
        \left\|
            \frac{1}{h_n}Q_nM_n^{-1}J
        \right\|_{\rm F}^2
        \mathbf 1_{\mathcal G_n}
    \right]
    \leq
    C\left(h_n^4+\frac{h_n^2}{n}\right)
    \leq
    C\left(h_n^4+\frac{1}{nh_n^2}\right),
\end{equation}
where the last inequality follows from $h_n^2\leq \bar h^4/h_n^2$, with $\bar h^4$ absorbed into $C$.
We now establish the bounds collected in~\eqref{eq:adaptive_remainder_ingredients}.
Partition $Q$ as $Q=\begin{pmatrix}Q^{(0)}&Q^{(1)}\end{pmatrix}$, where
\[
    Q^{(0)}\triangleq
    \int_{\mathcal S}K(u)\wt g(u)R_z(u)\,{\rm d}u,
    \qquad
    Q^{(1)}\triangleq
    \int_{\mathcal S}K(u)\wt g(u)R_z(u)u^\top\,{\rm d}u.
\]
The uniform bounds on the second- and third-order derivatives of $c$ imply that there exist constants $L_2,L_3>0$ such that, for all $z\in\mathcal X$ and $u\in\mathcal S$,
\begin{equation}\label{eq:adaptive_taylor_expansion}
    R_z(u)
    =
    \frac{h_n^2}{2}\nabla^2c(z)[u,u]+\rho_z(u),
    \qquad
    \|\nabla^2c(z)[u,u]\|\leq L_2\|u\|^2,
    \qquad
    \|\rho_z(u)\|\leq L_3h_n^3\|u\|^3.
\end{equation}
Substituting~\eqref{eq:adaptive_taylor_expansion} into $Q^{(0)}$ gives
\begin{equation}\label{eq:adaptive_Q_intercept}
    \|Q^{(0)}\|_{\rm F}
    \leq
    \frac{L_2h_n^2}{2}
    \int_{\mathcal S}K(u)\wt g(u)\|u\|^2\,{\rm d}u
    +L_3h_n^3
    \int_{\mathcal S}K(u)\wt g(u)\|u\|^3\,{\rm d}u\leq
    \frac{L_2K_{\max}R_{\mathcal S}^2h_n^2}{2}
    +L_3K_{\max}R_{\mathcal S}^3h_n^3.
\end{equation}
For $Q^{(1)}$, \eqref{eq:adaptive_taylor_expansion} gives
\[
\begin{aligned}
    Q^{(1)}
    ={}&
    \frac{h_n^2}{2}
    \int_{\mathcal S}
    K(u)\wt g(u)\nabla^2c(z)[u,u]u^\top\,{\rm d}u +
    \int_{\mathcal S}
    K(u)\wt g(u)\rho_z(u)u^\top\,{\rm d}u.
\end{aligned}
\]
Since $K$, $\wt g$, and $\nabla^2c(z)[u,u]$ are even in $u$, the first integral vanishes.
Hence,
\begin{equation}\label{eq:adaptive_Q_slope}
    \|Q^{(1)}\|_{\rm F}
    \leq
    L_3h_n^3
    \int_{\mathcal S}K(u)\wt g(u)\|u\|^4\,{\rm d}u
    \leq
    L_3K_{\max}R_{\mathcal S}^4h_n^3.
\end{equation}
Define
$
    B_0\triangleq
    \frac{L_2K_{\max}R_{\mathcal S}^2}{2}
    +L_3\bar hK_{\max}R_{\mathcal S}^3$ and $
    B_1\triangleq
    L_3K_{\max}R_{\mathcal S}^4
$.
Since $h_n\leq\bar h$,~\eqref{eq:adaptive_Q_intercept} and~\eqref{eq:adaptive_Q_slope} imply
\begin{equation}\label{eq:adaptive_Q_bound}
    \|Q\|_{\rm F}^2
    =
    \|Q^{(0)}\|_{\rm F}^2+\|Q^{(1)}\|_{\rm F}^2
    \leq
    \left(B_0^2+B_1^2\bar h^2\right)h_n^4.
\end{equation}
The symmetry of $K$ and $\wt g$ also implies
\[
    M=
    \begin{pmatrix}
        \mu&\mathbf 0_d^\top\\
        \mathbf 0_d&\Sigma
    \end{pmatrix},
    \qquad
    M^{-1}J=
    \begin{pmatrix}
        \mathbf 0_d^\top\\
        \Sigma^{-1}
    \end{pmatrix}.
\]
Since $M\succeq\kappa\mathbb I_{d+1}$, we have $\Sigma\succeq\kappa\mathbb I_d$ and $\|\Sigma^{-1}\|\leq1/\kappa$.
It follows from~\eqref{eq:adaptive_Q_slope} that
\begin{equation}\label{eq:adaptive_population_bias}
    \|QM^{-1}J\|_{\rm F}
    =
    \|Q^{(1)}\Sigma^{-1}\|_{\rm F}
    \leq
    \frac{B_1}{\kappa}h_n^3.
\end{equation}
We next bound $\mathbb E\left[\|Q_n-Q\|_{\rm F}^2\right]$.
Define
$
    Y_i\triangleq
    K(U_i)R_z(U_i)\psi(U_i)^\top$ and $
    Z_i\triangleq Y_i-Q
$.
Then $Q_n-Q=n^{-1}\sum_{i=1}^nZ_i$, and the matrices $Z_i$ are independent with zero mean.
By~\eqref{eq:adaptive_taylor_expansion},
    $\|R_z(U_i)\|
    \leq
    \left(
        \frac{L_2R_{\mathcal S}^2}{2}
        +L_3\bar hR_{\mathcal S}^3
    \right)h_n^2$.
Since $\|\psi(U_i)\|\leq\sqrt{1+R_{\mathcal S}^2}$, define
    $B_Y\triangleq
    K_{\max}\sqrt{1+R_{\mathcal S}^2}
    \left(
        \frac{L_2R_{\mathcal S}^2}{2}
        +L_3\bar hR_{\mathcal S}^3
    \right)$.
Then $\|Y_i\|_{\rm F}\leq B_Yh_n^2$ almost surely.
Using the Frobenius inner product and the independence and zero means of the $Z_i$, we obtain
\begin{equation}\label{eq:adaptive_Q_fluctuation}
\begin{aligned}
    \mathbb E\left[\|Q_n-Q\|_{\rm F}^2\right]
    &=
    \frac{1}{n}
    \mathbb E\left[\|Z_1\|_{\rm F}^2\right]=
    \frac{1}{n}\left(\mathbb E\left[\|Y_1\|_{\rm F}^2\right]-\|Q\|_{\rm F}^2\right)\leq
    \frac{1}{n}
    \mathbb E\left[\|Y_1\|_{\rm F}^2\right]
    \leq
    \frac{B_Y^2h_n^4}{n}.
\end{aligned}
\end{equation}
Finally, define
$
    W_i\triangleq
    K(U_i)\psi(U_i)\psi(U_i)^\top$ and $
    \Delta_i\triangleq W_i-M
$.
Then $M_n-M=n^{-1}\sum_{i=1}^n\Delta_i$, where the matrices $\Delta_i$ are independent with zero mean.
Moreover,
    $\|W_i\|_{\rm F}
    =
    K(U_i)\|\psi(U_i)\|^2
    \leq
    K_{\max}\left(1+R_{\mathcal S}^2\right)$.
Hence,
\begin{equation}\label{eq:adaptive_M_fluctuation}
\begin{aligned}
    \mathbb E\left[\|M_n-M\|^2\right]
    &\leq
    \mathbb E\left[\|M_n-M\|_{\rm F}^2\right]=
    \frac{1}{n}
    \mathbb E\left[\|\Delta_1\|_{\rm F}^2\right]=
    \frac{1}{n}\left(\mathbb E\left[\|W_1\|_{\rm F}^2\right]-\|M\|_{\rm F}^2\right)\leq
    \frac{K_{\max}^2\left(1+R_{\mathcal S}^2\right)^2}{n}.
\end{aligned}
\end{equation}
On $\mathcal G_n$, the resolvent identity and~\eqref{eq:adaptive_inverse_bound} give
\begin{equation}\label{eq:adaptive_inverse_perturbation}
    \|M_n^{-1}-M^{-1}\| = \|M_n^{-1}(M-M_n)M^{-1}\| \leq
    \frac{2}{\kappa^2}\|M_n-M\|.
\end{equation}
Equations~\eqref{eq:adaptive_Q_bound},~\eqref{eq:adaptive_population_bias},~\eqref{eq:adaptive_Q_fluctuation},~\eqref{eq:adaptive_M_fluctuation}, and~\eqref{eq:adaptive_inverse_perturbation} establish all the bounds used in~\eqref{eq:adaptive_remainder_ingredients}.

Since the operator norm is bounded by the Frobenius norm, the inequality $\|A+B\|^2\leq2\|A\|^2+2\|B\|^2$, together with~\eqref{eq:adaptive_gradient_decomposition},~\eqref{eq:adaptive_noise_bound}, and~\eqref{eq:adaptive_remainder_bound}, yields
\begin{equation}\label{eq:adaptive_good_mse}
    \mathbb E\left[
        \|G_n(z)-\nabla c(z)\|^2
        \mathbf 1_{\mathcal G_n}
    \right]
    \leq
    C\left(h_n^4+\frac{1}{nh_n^2}\right).
\end{equation}
Every constant used above is independent of $z$, and hence the bound holds uniformly over $z\in\mathcal X$.
It remains to control the complement of the good event $\mathcal{G}_n$.
Lemma~\ref{lem:truncated_gradient_estimator} gives
    $\|\wh\nabla c_n(z)-\nabla c(z)\|
    \leq
    2\|G_n(z)-\nabla c(z)\|$,
while Assumption~\ref{assumpt:SP-DD}(d) and the definition of the estimator $\wh\nabla c_n(z)$ give
    $\|\wh\nabla c_n(z)-\nabla c(z)\|
    \leq
    \wh{\ell}_c+\ell_c$.
Therefore, by~\eqref{eq:adaptive_good_event} and~\eqref{eq:adaptive_good_mse},
\[
\begin{aligned}
    \mathbb E\left[
        \|\wh\nabla c_n(z)-\nabla c(z)\|^2
    \right]
    &\leq
    4\mathbb E\left[
        \|G_n(z)-\nabla c(z)\|^2
        \mathbf 1_{\mathcal G_n}
    \right]
    +(\wh{\ell}_c+\ell_c)^2\mathbb P(\mathcal G_n^c)\\
    &\leq
    C_0\left(h_n^4+\frac{1}{nh_n^2}\right)
    +C_1\exp(-c_0n).
\end{aligned}
\]
Taking the supremum over $z\in\mathcal X$ proves~\eqref{eq:adaptive_uniform_mse}.

\subsection{Proof of Proposition~\ref{prop:opt_error}}
\label{sec:proof_opt_error}

Let $\eta_0^t$ be a sample generated from $\mathbb P_\xi(z^t)$ independently of $W_{m_t}^t$ conditional on $B^0,\ldots,B^{t-1}$ and $S_{n_t}^t$. For notational convenience, define
\[
    \wh F^{n_t}(x,\eta;z^t) \triangleq \varphi\left(x,\eta+\wh{\nabla}c_{n_t}^t(z^t)(x-z^t)\right).
\]
With $\mathbb{E}_{t,0}[\,\cdot\,]$ denoting the expectation with respect to $\eta_0^t$ conditional on all other random variables, we have
\begin{equation*}
\begin{aligned}
    &\mathrel{\phantom{=}}
    \left|\mathbb E_{B^t}\left[
        \wh F^{n_t,m_t}(z^{t+1},W_{m_t}^t;z^t)-f(z^{t+1})
    \right]\right| \\
    &\leq
    \left|\mathbb E_{B^t}\left[
        \wh F^{n_t,m_t}(z^{t+1},W_{m_t}^t;z^t)
        -\mathbb E_{t,0}\left[\wh F^{n_t}(z^{t+1},\eta_0^t;z^t)\right]
    \right]\right| + \left|\mathbb E_{B^t}\mathbb E_{t,0}\left[
        \wh F^{n_t}(z^{t+1},\eta_0^t;z^t)-f(z^{t+1})
    \right]\right|.
\end{aligned}
\end{equation*}
Noticing that the function $\wh F^{n_t}(\,\cdot\,,\eta;z^t)$ is $\ell_\varphi(1+\wh{\ell}_c)$-Lipschitz continuous and $\lambda(1+\wh{\ell}_c)^2$-weakly convex for every $\eta\in\mathbb R^\ell$, and $\eta_0^t$ is independently and identically distributed with the samples in $W^t_{m_t}$ conditional on $B^0, \dots, B^{t-1}, S^t_{n_t}$, we apply \citet[Theorem 3.2]{deng2021minibatch} to obtain
\begin{equation*}
\begin{aligned}
    &\left|\mathbb E_{B^t}\left[
        \wh F^{n_t,m_t}(z^{t+1},W_{m_t}^t;z^t)
        -\mathbb E_{t,0}\left[\wh F^{n_t}(z^{t+1},\eta_0^t;z^t)\right]
    \right]\right| \leq
    \frac{2\ell_\varphi^2(1+\wh{\ell}_c)^2}
    {m_t(\alpha_t-\lambda(1+\wh{\ell}_c)^2)}.
\end{aligned}
\end{equation*}
For the second term, we have
\begin{equation*}
\begin{aligned}
    &\mathrel{\phantom{=}}
    \left|\mathbb E_{B^t}\mathbb E_{t,0}\left[
        \wh F^{n_t}(z^{t+1},\eta_0^t;z^t)-f(z^{t+1})
    \right]\right| \\
    &\leq
    \left|\mathbb E_{B^t}\mathbb E_{t,0}\left[
        \wh F^{n_t}(z^{t+1},\eta_0^t;z^t)
        -\wt F(z^{t+1},\eta_0^t;z^t)
    \right]\right| +
    \left|\mathbb E_{B^t}\mathbb E_{t,0}\left[
        \wt F(z^{t+1},\eta_0^t;z^t)-f(z^{t+1})
    \right]\right| \\
    &\leq
    \ell_\varphi\mathbb E_{B^t}\left[
        \left\|\wh{\nabla}c_{n_t}^t(z^t)-\nabla c(z^t)\right\|
        \left\|z^{t+1}-z^t\right\|
    \right]
    +\frac{\ell_\varphi L_c}{2}
    \mathbb E_{B^t}\left[\left\|z^{t+1}-z^t\right\|^2\right] \\
    &\leq
    \ell_\varphi e_1(n_t)
    \sqrt{\mathbb E_{B^t}\left[\left\|z^{t+1}-z^t\right\|^2\right]}
    +\frac{\ell_\varphi L_c}{2}
    \mathbb E_{B^t}\left[\left\|z^{t+1}-z^t\right\|^2\right],
\end{aligned}
\end{equation*}
where the second inequality follows from the Lipschitz continuity of $\varphi(\,\cdot\,,\,\cdot\,)$ and Proposition~\ref{prop:pro_theo_surro}(b), and the third inequality follows from the Cauchy-Schwarz inequality. Combining the above inequalities yields the claimed result.

\subsection{Proof of Proposition~\ref{prop:des_property}}
\label{sec:proof_des_property}

By the strong convexity of $f(\,\cdot\,)+\frac{\rho}{2}\|\,\cdot\, - z^t\|^2$,
$\wh{F}^{n_t, m_t} \left( \,\cdot\,, W_{m_t}^t; z^t \right) + \frac{\alpha_t}{2}\left\|\,\cdot\, - z^t\right\|^2$ and definitions of $\wh{z}^{\,t}$, $z^{t+1}$,
we have
\begin{equation*}
\begin{aligned}
    \wh{F}^{n_t,m_t} \left( z^{t+1}, W_{m_t}^t; z^t \right) + \frac{\alpha_t}{2}\left\|z^{t+1} - z^t\right\|^2 &\leq \wh{F}^{n_t,m_t} \left( \wh{z}^{\,t}, W_{m_t}^t; z^t \right) + \frac{\alpha_t}{2}\left\|\wh{z}^{\,t} - z^t\right\|^2 - \frac{\alpha_t - \wh{\lambda}}{2} \left\|\, \wh{z}^{\,t}-z^{t+1} \,\right\|^2, \\
    f\left( \wh{z}^{\,t} \right) + \frac{\rho}{2}\left\|\, \wh{z}^{\,t}-z^t \,\right\|^2 &\leq f\left( z^{t+1} \right) + \frac{\rho}{2}\left\|\, z^{t+1}-z^t \,\right\|^2 - \frac{\rho - \ell_\varphi L_c - \wt{\lambda}}{2}\left\|\, z^{t+1}-\wh{z}^{\,t} \,\right\|^2.
\end{aligned}
\end{equation*}
Adding the above inequalities and rearranging terms, with $\tau = \ell_\varphi L_c + \wh{\lambda} + \wt{\lambda}$, we deduce
\begin{equation*}
\begin{aligned}
    &\mathrel{\phantom{=}} \frac{\alpha_t + \rho - \tau}{2} \left\|\, \wh{z}^{\,t}-z^{t+1} \,\right\|^2 + \frac{\rho - \alpha_t}{2}\left\|\, \wh{z}^{\,t}-z^t \,\right\|^2 + \frac{\alpha_t - \rho}{2} \left\|\, z^{t+1}-z^t \,\right\|^2\\
    &\leq \wh{F}^{n_t,m_t} \left( \wh{z}^{\,t}, W_{m_t}^t; z^t \right) - f\left( \wh{z}^{\,t} \right) + f\left( z^{t+1} \right) - \wh{F}^{n_t,m_t} \left( z^{t+1}, W_{m_t}^t; z^t \right) .
\end{aligned}
\end{equation*}
By Propositions~\ref{prop:cvx_prac_surro} and~\ref{prop:opt_error},
taking expectations conditioned on $B^0, \ldots, B^{t-1}$ on both sides of the above inequality yields
\begin{equation*}
\begin{aligned}
    &\mathrel{\phantom{=}} \frac{\alpha_t + \rho - \tau}{2} \mathbb{E}_{B^t} \left[\left\|\, \wh{z}^{\,t}-z^{t+1} \,\right\|^2\right] + \frac{\rho - \alpha_t}{2}\left\|\, \wh{z}^{\,t}-z^t \,\right\|^2 + \frac{\alpha_t - \rho}{2} \mathbb{E}_{B^t} \left[\left\|\, z^{t+1}-z^t \,\right\|^2\right] \\
    &\leq \ell_\varphi e_1(n_t) \left\|\, \wh{z}^{\,t}-z^t \,\right\| + \frac{\ell_\varphi L_c}{2} \left\|\, \wh{z}^{\,t}-z^t \,\right\|^2 + \frac{2\ell_\varphi^2 (1 + \wh{\ell}_c)^2}{m_t(\alpha_t - \wh{\lambda})} + \ell_\varphi e_1(n_t) \sqrt{\mathbb{E}_{B^t}\left[\left\| z^{t+1} - z^t \right\|^2 \right]} \\
    &\epc + \frac{\ell_\varphi L_c}{2} \mathbb{E}_{B^t}\left[\, \left\|\, z^{t+1}-z^t \,\right\|^2 \,\right] .
\end{aligned}
\end{equation*}
Rearranging terms on the above inequality, we have
\begin{equation*}
\begin{aligned}
    &\mathrel{\phantom{=}} \frac{2\rho - \ell_\varphi L_c - \tau}{2} \left\|\, \wh{z}^{\,t} - z^t \,\right\|^2 - \ell_\varphi e_1(n_t) \left\|\, \wh{z}^{\,t} - z^t \,\right\| \\
    &\leq \frac{\alpha_t + \rho - \tau}{2} \left(\left\|\, \wh{z}^{\,t}-z^t \,\right\|^2 - \mathbb{E}_{B^t}\left[\, \left\|\, z^{t+1}-\wh{z}^{\,t} \,\right\|^2 \,\right]\right) + \frac{2\ell_\varphi^2 (1 + \wh{\ell}_c)^2}{m_t(\alpha_t - \wh{\lambda})} \\
    & \epc - \frac{\alpha_t - \rho - \ell_\varphi L_c}{2} \mathbb{E}_{B^t}\left[\left\| z^{t+1} - z^t \right\|^2 \right] + \ell_\varphi e_1(n_t) \sqrt{\mathbb{E}_{B^t}\left[\left\| z^{t+1} - z^t \right\|^2 \right]} \\
    &\leq \frac{\alpha_t + \rho - \tau}{2} \left(\left\|\, \wh{z}^{\,t}-z^t \,\right\|^2 - \mathbb{E}_{B^t}\left[\, \left\|\, z^{t+1}-\wh{z}^{\,t} \,\right\|^2 \,\right]\right) + \frac{2\ell_\varphi^2 (1 + \wh{\ell}_c)^2}{m_t(\alpha_t - \wh{\lambda})} + \frac{\ell_\varphi^2 e_1^2(n_t)}{2(\alpha_t - \rho - \ell_\varphi L_c)} ,
\end{aligned}
\end{equation*}
where the second inequality follows by $\alpha_t > \rho + \ell_\varphi L_c$ and taking maximum of the function $x \mapsto -\frac{\alpha_t - \rho - \ell_\varphi L_c}{2} x^2 + \ell_\varphi e_1(n_t) x$ over $x \geq 0$.
Multiplying $\rho / (\alpha_t + \rho - \tau)$ on both sides of the above inequality, we obtain
{\small
\begin{equation*}
\begin{aligned}
    &\mathrel{\phantom{=}} \frac{\rho(2\rho - \ell_\varphi L_c - \tau)}{2(\alpha_t + \rho - \tau)} \left\|\, \wh{z}^{\,t} - z^t \,\right\|^2 - \frac{\rho \ell_\varphi e_1(n_t)}{\alpha_t + \rho - \tau} \left\|\, \wh{z}^{\,t} - z^t \,\right\| \\
    &\leq \frac{\rho}{2} \left(\left\| \wh{z}^{\,t}-z^t \right\|^2 - \mathbb{E}_{B^t}\left[ \left\|\, z^{t+1}-\wh{z}^{\,t} \,\right\|^2 \right]\right) + \frac{\rho}{\alpha_t + \rho - \tau} \left[\frac{2 \ell_\varphi^2 (1 + \wh{\ell}_c)^2}{m_t(\alpha_t - \wh{\lambda})} + \frac{\ell_\varphi^2 e_1^2(n_t)}{2(\alpha_t - \rho - \ell_\varphi L_c)}\right] \\
    &= f\left( \wh{z}^{\,t} \right) + \frac{\rho}{2} \left\|\, \wh{z}^{\,t}-z^t \,\right\|^2 - f\left( \wh{z}^{\,t} \right) - \frac{\rho}{2} \mathbb{E}_{B^t}\left[\, \left\|\, z^{t+1}-\wh{z}^{\,t} \,\right\|^2 \,\right]  + \frac{\rho}{\alpha_t + \rho - \tau} \left[\frac{2 \ell_\varphi^2 (1 + \wh{\ell}_c)^2}{m_t(\alpha_t - \wh{\lambda})} + \frac{\ell_\varphi^2 e_1^2(n_t)}{2(\alpha_t - \rho - \ell_\varphi L_c)}\right] \\
    &\leq f_{1/\rho} \left(z^t\right) - \mathbb{E}_{B^t} \left[\, f_{1/\rho} \left(z^{t+1}\right) \,\right] + \frac{\rho}{\alpha_t + \rho - \tau} \left[\frac{2 \ell_\varphi^2 (1 + \wh{\ell}_c)^2}{m_t(\alpha_t - \wh{\lambda})} + \frac{\ell_\varphi^2 e_1^2(n_t)}{2(\alpha_t - \rho - \ell_\varphi L_c)}\right] ,
\end{aligned}
\end{equation*}
}
where the last inequality is obtained from the definition of Moreau envelope. We then obtain the claimed inequality~\eqref{eq:des_property} by taking expectations with respect to $B^0$, $\dots$, $B^{t-1}$ on both sides of the above inequality.

\subsection{Proof of Theorem~\ref{thm:nonasym_basic}}
\label{sec:proof_nonasym_basic}

By adding the inequality \eqref{eq:des_property} from $t=0$ to $t=T$ and noticing that $\mathbb{E} \left[\, f_{1/\rho} \left(z^{T+1}\right) \,\right] \geq \min_{z\in\mathcal{X}} f(z)$, we obtain
\begin{equation*}
\begin{aligned}
    &\sum_{t=0}^T \frac{\rho}{\alpha_t + \rho - \tau} \mathbb{E}\left[\frac{2\rho - \ell_\varphi L_c - \tau}{2} \left\|\, \wh{z}^{\,t} - z^t \,\right\|^2 - \ell_\varphi e_1(n_t) \left\|\, \wh{z}^{\,t} - z^t \,\right\|\right] \\
    &\qquad\leq \Delta_0 + \sum_{t=0}^T \frac{\rho}{\alpha_t + \rho - \tau} \left[\frac{2 \ell_\varphi^2 (1 + \wh{\ell}_c)^2}{m_t(\alpha_t - \wh{\lambda})} + \frac{\ell_\varphi^2 e_1^2(n_t)}{2(\alpha_t - \rho - \ell_\varphi L_c)}\right] .
\end{aligned}
\end{equation*}
By definition of the random variable $t^*$, dividing by $\sum_{t=0}^T \frac{\rho}{\alpha_t + \rho - \tau}$ on both sides of the above inequality yields
{\small
\begin{equation*}
\begin{aligned}
    &\mathbb{E}\left[\frac{2\rho - \ell_\varphi L_c - \tau}{2} \left\|\, \wh{z}^{\,t^*} - z^{t^*} \,\right\|^2 - \ell_\varphi e_1(n_{t^*}) \left\|\, \wh{z}^{\,t^*} - z^{t^*} \,\right\|\right] \leq \frac{\Delta_0 + \sum_{t=0}^T \frac{\rho}{\alpha_t + \rho - \tau} \left[\frac{2 \ell_\varphi^2 (1 + \wh{\ell}_c)^2}{m_t(\alpha_t - \wh{\lambda})} + \frac{\ell_\varphi^2 e_1^2(n_t)}{2(\alpha_t - \rho - \ell_\varphi L_c)}\right]}{\sum_{t=0}^T \frac{\rho}{\alpha_t + \rho - \tau}}.
\end{aligned}
\end{equation*}
}
Denote by $\Delta(T)$ the right-hand side of the above inequality. We deduce
\begin{equation*}
\begin{aligned}
    \Delta(T) &\geq \mathbb{E}\left[\frac{2\rho - \ell_\varphi L_c - \tau}{2} \left\|\, \wh{z}^{\,t^*} - z^{t^*} \,\right\|^2 - \ell_\varphi e_1(n_{t^*}) \left\|\, \wh{z}^{\,t^*} - z^{t^*} \,\right\|\right] \\
    &= \mathbb{E} \left[\frac{2\rho - \ell_\varphi L_c - \tau}{2} \left(\left\|\, \wh{z}^{\,t^*} - z^{t^*} \,\right\| - \frac{\ell_\varphi e_1(n_{t^*})}{2\rho - \ell_\varphi L_c - \tau}\right)^2 - \frac{\ell_\varphi^2 e_1^2 (n_{t^*})}{2(2\rho - \ell_\varphi L_c - \tau)}\right] \\
    &\geq \frac{2\rho - \ell_\varphi L_c - \tau}{2} \left(\mathbb{E}\left[\left\|\, \wh{z}^{\,t^*} - z^{t^*} \,\right\| - \frac{\ell_\varphi e_1(n_{t^*})}{2\rho - \ell_\varphi L_c - \tau}\right]\right)^2 - \frac{\ell_\varphi^2 \mathbb{E}[e_1^2 (n_{t^*})]}{2(2\rho - \ell_\varphi L_c - \tau)} .
\end{aligned}
\end{equation*}
Rearranging terms, taking square roots and noticing that $\mathbb{E}[e_1(n_{t^*})] \leq \sqrt{\mathbb{E}\left[e^2_1(n_{t^*})\right]}$ yield
\begin{equation*}
\begin{aligned}
    \mathbb{E}\left\|\, \wh{z}^{\,t^*} - z^{t^*} \,\right\| \leq \frac{\ell_\varphi \sqrt{\mathbb{E}\left[e^2_1(n_{t^*})\right]}}{2\rho - \ell_\varphi L_c - \tau} + \sqrt{\frac{2\Delta(T)}{2\rho - \ell_\varphi L_c - \tau} + \frac{\ell_\varphi^2 \mathbb{E}[e_1^2(n_{t^*})]}{(2\rho - \ell_\varphi L_c - \tau)^2}}.
\end{aligned}
\end{equation*}
Denote $a \triangleq \frac{2\Delta(T)}{2\rho - \ell_\varphi L_c - \tau}$ and $b \triangleq \frac{\ell_\varphi^2 \mathbb{E}[e_1^2(n_{t^*})]}{(2\rho - \ell_\varphi L_c - \tau)^2}$. The inequality $2\sqrt{b(a+b)}\leq b+(a+b)$ gives $\left(\sqrt{b} + \sqrt{a + b}\right)^2 \leq 2a + 4b$. Thus we obtain
\begin{equation*}
\begin{aligned}
    \mathbb{E}\left\|\, \wh{z}^{\,t^*} - z^{t^*} \,\right\| \leq \sqrt{\frac{4\Delta(T)}{2\rho - \ell_\varphi L_c - \tau} + \frac{4\ell_\varphi^2 \mathbb{E}[e_1^2(n_{t^*})]}{(2\rho - \ell_\varphi L_c - \tau)^2}}.
\end{aligned}
\end{equation*}
Unfolding the definitions of $\Delta(T)$ and $\mathbb{E}[e_1^2(n_{t^*})]$ yields the claimed result.

\subsection{Proof of Corollary~\ref{coroll:vari}}
\label{sec:proof_coroll_vari}

By Theorem~\ref{thm:nonasym_basic} and Assumption~\ref{assumpt:estimators}, regarding the algorithm parameters $\{\alpha_t\}$, $\{m_t\}$ and $\{n_t\}$, we have
\begin{equation*}
\begin{aligned}
    \mathbb{E}\left\|\, \wh{z}^{\,t^*} - z^{t^*} \,\right\| = O\left(\sqrt{\frac{1 + \sum_{t=0}^T\left(\frac{1}{m_t\alpha_t^2} + \frac{1}{n_t^{2r}\alpha_t} + \frac{1}{n_t^{2r}\alpha_t^2}\right)}{\sum_{t=0}^T \frac{1}{\alpha_t}}}\right).
\end{aligned}
\end{equation*}
Under the parameter settings $\alpha_t = \alpha_0 (t+1)^b$, $m_t = \lceil m_0(t+1)^j \rceil$, $n_t = \lceil n_0(t+1)^k \rceil$, denote $q \triangleq \min\{j + 2b, 2kr + b\} \geq b$. We have
\begin{equation*}
\begin{aligned}
    \mathbb{E}\left\|\, \wh{z}^{\,t^*} - z^{t^*} \,\right\| = \begin{cases}
        O\left( \sqrt{\frac{1}{T^{q - b}}} \right), &\text{ if } q < 1, \\
        \wt{O}\left( \sqrt{\frac{1}{T^{1-b}}} \right), &\text{ if } q \geq 1.
    \end{cases}
\end{aligned}
\end{equation*}
Note that $q \geq 1$ is equivalent to $j \geq j^*(b)$ and $k \geq k^*(b)$. The setting with $j \geq j^*(b)$ and $k \geq k^*(b)$ can not yield better sample complexity than $j = j^*(b)$ and $k = k^*(b)$ since the iteration complexities of these two settings are the same. When $j = j^*(b)$ and $k=k^*(b)$, by the above result, the iteration complexity is $T=\wt{O}\left(\wt{\nu}^{-\frac{2}{1-b}}\right)$, which yields the sample complexity of
\begin{equation*}
\begin{aligned}
    \sum_{t=0}^T(m_t + n_t) = O\left(T^{1+j^*(b)} + T^{1+k^*(b)}\right) = \wt{O}\left(\wt{\nu}^{-\frac{2(1+j^*(b))}{1-b}} + \wt{\nu}^{-\frac{2(1+k^*(b))}{1-b}}\right) = \wt{O}\left(\wt{\nu}^{-4} + \wt{\nu}^{-\frac{2}{1-b} - \frac{1}{r}}\right).
\end{aligned}
\end{equation*}
We have $4 < \frac{2}{1-b} + \frac{1}{r}$ since $r < \frac{1}{2}$.
When $q < 1$, the iteration complexity is $T=O\left(\wt{\nu}^{-\frac{2}{q-b}}\right)$, which yields the sample complexity of
\begin{equation*}
\begin{aligned}
    \sum_{t=0}^T(m_t + n_t) = O\left(T^{1+j} + T^{1+k}\right) = O\left(\wt{\nu}^{-\frac{2(1+j)}{q-b}} + \wt{\nu}^{-\frac{2(1+k)}{q-b}}\right).
\end{aligned}
\end{equation*}
In this case, if $k \geq k^*(b)$, we have $\frac{2(1+k)}{q-b} \geq \frac{2(1+k^*(b))}{q-b} > \frac{2(1+k^*(b))}{1-b}$. If $k \in [0, k^*(b))$, we have
\begin{equation*}
\begin{aligned}
    \frac{2(1+k)}{q-b} - \frac{2(1+k^*(b))}{1-b} &= \frac{2(1+k)(1-b) - 2(1+k^*(b))(q-b)}{(q-b)(1-b)} \geq \frac{2(1+k)(1-b) - 4kr(1+k^*(b))}{(q-b)(1-b)} \\
    &= \frac{2(1-b+k-kb-2kr-k+kb)}{(q-b)(1-b)} = \frac{2(1-b-2kr)}{(q-b)(1-b)} > 0,
\end{aligned}
\end{equation*}
where the first inequality is obtained by $q \leq 2kr + b$, and the second inequality is obtained by $k < k^*(b)$. Similarly, we obtain that $\frac{2(1+j)}{q-b} \geq \frac{2(1+j^*(b))}{1-b}$ for all $j \geq 0$. Therefore, the sample complexity of the parameter setting with $q < 1$ can not be better than that of the setting with $j=j^*(b)$ and $k=k^*(b)$, and we complete the proof.

\section{Verification of assumptions on statistical estimators}
\label{sec:estimation}

This appendix supplements the statistical estimation discussion in Section~\ref{sec:static_design}.
We first show by the following lemma that any gradient estimator satisfying Assumption~\ref{assumpt:estimators} can be truncated to satisfy Assumption~\ref{assu:bnd_gra_est} without changing its convergence rate.
We then examine the attainable values of $r$ in Assumption~\ref{assumpt:estimators} for several nonparametric gradient estimators under static random designs other than the LLR estimator discussed in Section~\ref{sec:static_design}.

\begin{lemma}
\label{lem:truncated_gradient_estimator}
Suppose that Assumption~\ref{assumpt:SP-DD} holds.
Given a scalar $\wh{\ell}_c\geq \ell_c$, define
\[
    \wh{\nabla}c_n(z)
    \triangleq
    \min\left\{
        1,
        \frac{\wh{\ell}_c}{\|G_n(z)\|}
    \right\}
    G_n(z).
\]
Then for every $z\in\mathcal X$ and $n\geq1$, $\big\|\wh{\nabla}c_n(z)\big\| \leq \wh{\ell}_c$, and
\[
    \big\|\wh{\nabla}c_n(z)-\nabla c(z)\big\| \leq 2\big\|G_n(z)-\nabla c(z)\big\|.
\]
Consequently, if $G_n(z)$ satisfies Assumption~\ref{assumpt:estimators}, then the truncated estimate $\wh{\nabla}c_n(z)$ satisfies Assumptions~\ref{assumpt:estimators} and~\ref{assu:bnd_gra_est} with the same convergence rate $r$.
\end{lemma}

\begin{proof}
The boundedness of $\wh{\nabla}c_n(z)$ follows from its definition.
Since $c$ is differentiable and $\ell_c$-Lipschitz continuous
under Assumption~\ref{assumpt:SP-DD},
we have
    $\|\nabla c(z)\|
    \leq
    \ell_c
    \leq
    \wh{\ell}_c$ for all $z\in\mathcal X$.
If $\|G_n(z)\|\leq\wh{\ell}_c$, then
$\wh{\nabla}c_n(z)=G_n(z)$.
Otherwise,
\[
    \big\|
        \wh{\nabla}c_n(z)-G_n(z)
    \big\|
    =
    \|G_n(z)\|-\wh{\ell}_c \leq
    \|G_n(z)\|-\|\nabla c(z)\| \leq
    \big\|
        G_n(z)-\nabla c(z)
    \big\|,
\]
Therefore, in both cases,
\[
    \big\|
        \wh{\nabla}c_n(z)-\nabla c(z)
    \big\|
    \leq \big\|\wh{\nabla}c_n(z)-G_n(z)\big\| + \big\|G_n(z)-\nabla c(z)\big\| \leq
    2
    \big\|
        G_n(z)-\nabla c(z)
    \big\|,
\]
which completes the proof.
\end{proof}

We next discuss convergence rates for several nonparametric gradient estimation methods, including regression splines, $k$-nearest neighbors ($k$-NN), difference-based methods, and neural networks.
The result for regression splines directly verifies Assumption~\ref{assumpt:estimators}.
For the other methods, available convergence-rate results differ from Assumption~\ref{assumpt:estimators} in their error criteria or design settings.
Nevertheless, they can still provide useful information about the attainable rates.

Regression spline estimators are constructed from a set of spline functions through (weighted) least squares with or without penalty terms. Splines are piecewise polynomials joined at knots with prescribed continuity conditions.
For multivariate predictors, common extensions include tensor product splines and thin plate splines \citep{hastie2009elements}.
Suppose that all partial derivatives of the regression function of total order less than or equal to $p-1$ are Lipschitz continuous.
According to  \citet[Theorem 2.5]{ma2020nonparametric}, if $p > d/2$ and the number of interior knots per coordinate $K$ satisfies $K = \Theta(n^\beta)$ where $\beta\in(\frac{1}{2p+d}, \frac{1}{2d})$, then the gradient estimator based on tensor-product splines of order $p$ satisfies Assumption~\ref{assumpt:estimators} with $r = \frac{1}{2}(1 - (d+2)\beta)$. By choosing $\beta$ arbitrarily close to $1/(2p+d)$, the rate can be made arbitrarily close to $\frac{p-1}{2p+d}$.

The usual $k$-NN estimate averages the responses at the $k$ nearest design points \citep{biau_lectures_2015}.
To estimate the gradient, \citet{pmlr-v130-ausset21a} fit a penalized local linear approximation over these neighbors.
Their Theorem~1 gives a pointwise high-probability gradient error bound of order $n^{-1/(d+4)}$, up to logarithmic factors, when $k=\Theta\big(n^{4/(d+4)}\big)$.

Difference-based estimators estimate the gradient as a linear combination of responses at neighboring design points.
\citet{zhang2024difference} analyze such estimators for a $d$-dimensional equidistant tensor grid.
Their Theorems~2 and~4 yield a root MSE convergence rate $O(n^{-2/(d+6)})$ for interior and boundary grid points, respectively.
These rates apply only at grid points under a fixed equidistant design, whereas Assumption~\ref{assumpt:estimators} concerns uniform error over $\mathcal X$ under the static random design.

For neural networks, \citet{ding2025sobolev} analyze a ReQU estimator for a $C^s$-smooth regression function.
Their Theorem~V.4 gives an $O(n^{-s/(d+4s)}\log^2 n)$ rate for its expected squared gradient error averaged over the design distribution.

\bibliographystyle{informs2014}
\bibliography{ref/MyLibrary}

\end{document}